\documentclass[a4paper]{amsart}

\title[Dirichlet Stereohedra for Quarter Cubic Groups]
{On the Number of Facets of\\ Three-Dimensional Dirichlet Stereohedra IV: \\Quarter Cubic Groups}
\thanks{Research partially supported by  the Spanish Ministry of Education and Science, grant number
 MTM2008-04699-C03-02.}

\author{Pilar Sabariego
        \and
        Francisco Santos}
        \address{Departamento de Matem\'aticas, Estad\'{\i}stica y
        Computaci\'on, Universidad de Canta\-bria, 39005 Santander, Spain}
        \email{sabariego@gmail.com, francisco.santos@unican.es.}

\usepackage{latexsym}
\usepackage{amssymb}
\usepackage{amsmath}
\usepackage[dvips]{graphicx}

\usepackage{color}

\usepackage{multirow}

\newtheorem{theorem}{Theorem}[section]

\newtheorem{lemma}[theorem]{Lemma}
\newtheorem{definition}[theorem]{Definition}

\newcommand{\Z}{\mathbb Z}
\newcommand{\R}{\mathbb R}

\newcommand{\Isom}{\operatorname{Isom}}
\newcommand{\Vor}{\operatorname{Vor}}
\newcommand{\VorExt}{\operatorname{VorExt}}
\newcommand{\Infl}{\operatorname{Infl}}
\newcommand{\conv}{\operatorname{conv}}
\newcommand{\nor}{\operatorname{\mathcal N}}
\newcommand{\auxtes}{\mathcal T}

\begin{document}

\begin{abstract}

In this paper we finish the intensive study of three-dimensional Dirichlet stereohedra started by the second author and D.~Bochi\c{s}, who showed that they cannot have more than  80 facets, except perhaps for crystallographic space groups in the cubic system.

Taking advantage of the recent, simpler classification of three-dimensional
crystallographic groups by Conway, Delgado-Friedrichs, Huson and Thurston, in a previous paper we proved that Dirichlet stereohedra for any of the 27 ``full" cubic groups cannot have more than $25$ facets. 
Here we study the remaining ``quarter" cubic groups. With a computer-assisted method, our main result  is that Dirichlet stereohedra for the 8 quarter groups, hence for all three-dimensional crystallographic groups, cannot have more than $92$ facets.

\end{abstract}

\maketitle

\section{Introduction}

 This is the last in a series of four papers (see \cite{Bochis-Santos-2001, Bochis-Santos-2006, SabariegoSantos-full}) devoted to bounds on the number of facets that Dirichlet stereohedra in Euclidean 3-space can have.
 
 A {\em stereohedron} is any bounded convex polyhedron which tiles the space by the action of some crystallographic group. A {\em Dirichlet stereohedron} for a certain crystallographic group $G$ is the Voronoi region $\Vor_{Gp}(p)$ of a point $p\in \R^{3}$ in the Voronoi diagram of an orbit $Gp$. 

The study of the maximum number of facets for stereohedra is related to Hilbert's 18th problem, ``{\em Building up the space with congruent polyhedra"} (see \cite{Hilbert-1935, Milnor-1974}). Bie\-berbach (1910) and Reinhardt (1932) answered completly the first two of Hilbert's specific questions, but other problems related to {\em monohedral tessellations} (i.~e., tessellations whose tiles are congruent) remain open. An exhaustive account of this topic appeared in a survey article by Gr\"unbaum and Shephard  \cite{GruShe-art}, where our problem, \emph{to determine the maximum number of facets---or, at least, a ``good" upper bound---for Dirichlet stereohedra in $\R^{3}$}, is mentioned as an important one.
Previous results on this problem are:

\begin{itemize}
\item The { \em fundamental theorem of stereohedra} (Delaunay, 1961 \cite{Delone-1961}) asserts that a stereohedron of dimension $d$ for a crystallographic group $G$ with $a$ {\em aspects}  cannot have more than $2^{d}(a+1)-2$ facets. The number of aspects of a crystallographic group $G$ is the index of 
its translational subgroup. Delone's bound for three-dimensional groups, which have up to 48 aspects, is
390 facets.

\item The three-dimensional stereohedron with the maximum number of facets known so far was found in 1980 by P. Engel (see \cite{Engel} and \cite[p. 964]{GruShe-art}). It is a Dirichlet stereohedron with 38 facets, for the cubic group $I4_1 32$, with 24 aspects.

\end{itemize}

There is agreement among the experts (see \cite[page 214]{Engel}, \cite[page 960]{GruShe-art}, \cite[page 50]{SchSen-1997}) that Engel's sterohedron is much closer than Delone's upper bound to having the maximum possible number of facets. Our results confirm this.

In 2000, the second author and D. Bochi\c{s} gave upper bounds for the number of facets of Dirichlet stereohedra. They did this by dividing the 219 affine conjugacy classes of three-dimensional crystallographic groups into three blocks, and using different tools for each. Their main results are:

\begin{itemize}
\item Within the 100 crystallographic groups which contain reflection planes, the exact maximum number of facets is 18~%
\cite{Bochis-Santos-2001}.

\item Within the 97 non-cubic crystallographic groups without reflection planes, they found Dirichlet stereohedra with 32 facets and proved that no one can have more than 80. Moreover, they got upper bounds of 50 and 38 for all but, respectively, 9 and 21 of the groups \cite{Bochis-Santos-2006}.

\item They also considered cubic groups, but they were only able to prove an upper bound of 162 facets for them \cite{Bochis-1999}.
\end{itemize}

In \cite{SabariegoSantos-full} we improved the bound for 14 of the 22 cubic groups without reflections planes, the 14 ``full groups":

 \begin{theorem}
 \label{thm:main-full}
 Dirichlet stereohedra for full cubic groups cannot have more than
 25 facets. 
 \end{theorem}

In this paper we give an upper bound for the remaining cubic groups: the 8 ``quarter groups". It has to be noted that to get these bounds, contrary to the ones in the previous papers of this series, computers are used. The upper bound we obtain for each quarter group is shown in Table~\ref{table:main}. Columns (1) to (4) are the bounds obtained in different phases or our method, the column labeled ``Final" is our final bound. Globally, we get the following.

 \begin{theorem}
 \label{thm:main-quarter}
 Dirichlet stereohedra for quarter cubic groups cannot have more than
 92 facets. 
 \end{theorem}

 \begin{table}
\[
\begin{array}{|c|c|c|c|c|c|c|c|}
\hline
\text{$|G:Q|$} & \text{Aspects} & \text{Group} & \multicolumn{5}{c|}{\text{Our bounds}} \\
\hline    & & & \text{(1)} & \text{(2)}  & \text{(3)} & \text{(4)} &  \text{Final} \\
\hline 
8 & 48 & \nor(Q)=I\frac{4_{1}}{g}\overline{3}\frac{2}{d} & 519 & 155 & 100 & 68  & \textbf{68} \\
\hline 
\multirow{3}*{4}&\multirow{3}*{24}
& I4_{1}32 & 264 & 96  & 55 & &\textbf{55} \\
\cline{3-8} & & I\overline{4}3d & 257 & 78 &  & 76 & \textbf{76} \\
\cline{3-8} & & I\frac{2}{g}\overline{3} & 260 & 77 & & 57  &\textbf{57} \\
\hline 
\multirow{3}*{2} 
& 24 & P4_{1}32 & 135 &  & 92 & &\textbf{92} \\
\cline{3-8} &12 & I2'3 & 131 & 48 &  & 46 &\textbf{46} \\
\cline{3-8} & 24 &P\frac{2_{1}}{a}\overline{3} & 132 & & & 86 &\textbf{86} \\
\hline 1 & 12 & Q=P2_{1}3 & 69 & & &  &\textbf{69} \\
\hline
\end{array}
\]
\begin{enumerate}
\item[](1) Bounds after processing triad rotations.

\item[](2) Bounds after diad rotations with axes parallel to
the coordinate axes.

\item[](3) Bounds after diagonal diad rotations.

\item[](4) Bounds after intersecting with planar projections.

\end{enumerate}
\medskip
\caption{Bounds for the number of facets of Dirichlet stereohedra of quarter cubic groups}
\label{table:main}
\end{table}

For the sake of completeness, we include here
the full list of other crystallographic groups for which the bounds proved in this series of papers is bigger than 38 (Table~\ref{Table.numero-maximo-no-cubicos}). This list is the same as Table~2 in~\cite{Bochis-Santos-2006}.
\begin{table}[htb]
 \renewcommand{\arraystretch}{1.2}
\begin{center}
\begin{tabular}{|c|c|c|c|}
\hline
Group & Asp. & Bound \\

\hline
 $I\overline{4}c2$    & 8  & {\bf 40}  \\
 $P\frac{4_2}{n}\frac{2}{g}\frac{2}{c}$
                      & 16  & {\bf 40}  \\
 $R\overline{3}$      & 6   & {\bf 42} \\
 $R32$                & 6    & {\bf 42} \\
 $R3c$                & 6    & {\bf 42} \\
 $I4_1cd$             & 8  & {\bf 44}  \\
 $P3_12$              & 6   & {\bf 48} \\
 \hline
\end{tabular}
\ 
\begin{tabular}{|c|c|c|c|}
\hline
Group & Asp. & Bound \\

\hline
 $P3_112$             & 6   & {\bf 48} \\
 $P6_1$               & 6    & {\bf 48} \\
 $P4_122$             & 8   & {\bf 50}    \\
 $C\frac{2}{a}\frac{2}{c}\frac{2}{c}$
                      & 8   & {\bf 50} \\
 $I\frac{2}{a}\frac{2}{c}\frac{2}{c}$
                      & 8  & {\bf 50} \\
 $P4_12_12$           & 8   & {\bf 64} \\
 $I\frac{4_1}{g}$     & 8  & {\bf 70} \\
\hline
\end{tabular}
\ 
\begin{tabular}{|c|c|c|c|}
\hline
Group & Asp. & Bound \\

\hline
$I4_122$             & 8  & {\bf 70}   \\
 $I\overline{4}2d$    & 8  & {\bf 70}   \\
 $F\frac{2}{d}\frac{2}{d}\frac{2}{d}$
                      & 8  & {\bf 70}  \\
 $P 6_2 2 2$          & 12  & {\bf 78} \\
 $P 6_1 2 2$          & 12 & {\bf 78} \\
 $R \overline{3} \frac{2}{c}$
                      & 12  & {\bf 79} \\
 $I\frac{4_1}{g}\frac{2}{c}\frac{2}{d}$
                      & 16 & {\bf 80} \\
\hline
\end{tabular}

\end{center}
 \renewcommand{\arraystretch}{1}
\caption{Non-cubic groups where our upper bound is larger 
than 38}
\label{Table.numero-maximo-no-cubicos}
\end{table}

  \begin{theorem}
 \label{thm:main}
 Three dimensional Dirichlet stereohedra cannot have more than
 92 facets. They can possibly have more than 38 facets only in one of the 29 groups
 listed in tables~\ref{table:main} and~\ref{Table.numero-maximo-no-cubicos}.
 \end{theorem}

\section{Preliminaries and outline}

\subsection{``Full'' and ``quarter'' cubic groups}

Our division of cubic groups into ``full'' and ``quarter'' ones comes from the 
recent classification of three-dimensional crystallographic groups developed in \cite{Conway-etal-2001} by Conway et al. They divide crystallographic groups into ``reducible" and ``irreducible", were {\em irreducible groups} are those that do not have any invariant direction. It turns out that they coincide with the cubic groups of the classical classification. Conway et al.~define {\em odd subgroup} of an irreducible group $G$ as the one generated by the rotations of order three, and show that:

\begin{theorem}[Conway et al.~ \cite{Conway-etal-2001}]
 \label{thm:conway-etal}
 \begin{enumerate}
 \item There are only two possible odd subgroups of cubic groups, that we denote $F$ and $Q$.
 \item Both $F$ and $Q$ are normal in \ $\Isom(\R^3)$. Hence, every cubic group lies between its odd subgroup and the normalizer $\nor(F)$ and $\nor(Q)$ of it.
  \end{enumerate}
 \end{theorem}

 The second property reduces the enumeration of
 cubic space groups to the enumeration, up to conjugacy, of subgroups
 of the two finite groups $\nor(F)/F$ and $\nor(Q)/Q$. $\nor(Q)/Q$ is dihedral of order 8 and $\nor(F)/F$ has order 16 and contains a dihedral subgroup of index 2. 
 
The main difference between $F$ and $Q$ is that $Q$ only contains triad rotations whose axes are mutually disjoint, while some triad rotation axes in $F$ intersect one another. $Q$ is a subgroup of $F$ of order four and because of that Conway et al.~call {\em full groups} those with odd subgroup equal to $F$ and {\em quarter groups} those with odd subgroup equal to $Q$: Quarter groups contain only a quarter of the possible rotation axes. There are 27 full groups (14 of them without reflection planes) and 8 quarter groups (none with reflection planes, because $\nor(Q)$ does not contain reflections).

 \subsection{The structure of quarter cubic groups}
 \label{sec:quarter-groups}
 
  Throughout the paper we use the International Crystallographic Notation for three-dimensional crystallographic groups; see, e.~g., \cite{Lockwood-1978}.

All the quarter cubic groups have the same odd group $Q$, a group of type $P2_{1}3$. 
By definition, $Q$ is generated by triad rotations. More precisely, through each point $(x,y,z)\in (\Z/2)^{3}$ exactly one rotation axis passes, with vector:
 
 \[ 
 \begin {array}{l}
(1,1,1) \, \text{if} \,  x\equiv y \equiv z \, \text{(mod 1)} \\
(-1,1,1) \, \text{if} \, y\not \equiv x \equiv z \, \text{(mod 1)} \\
(1,-1,1) \, \text{if} \, z \not \equiv x \equiv y \, \text{(mod 1)} \\
(1,1,-1) \, \text{if} \, x \not \equiv y \equiv z \, \text{(mod 1)} \\
\end{array}
\]

\begin{figure}
    \begin{center}
      \includegraphics[width=8cm,height=7cm]{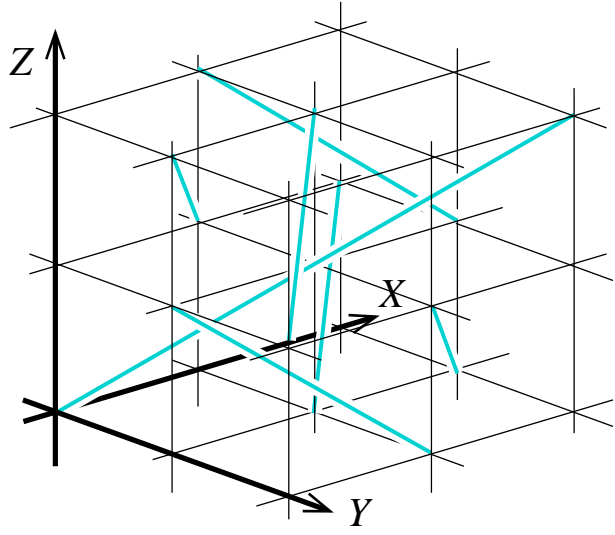}
    \end{center}
      \caption{The body-centered cubic lattice $I$. The triad rotations (in grey) together with translations of length two in the coordinate directions generate $Q$}
     \label{fig:red_giros}
\end{figure}

In other words, exactly one of the four diagonals of each primitive cubic cell of the lattice  $(\Z/2)^{3}$
is a rotation axis (see Fig~\ref{fig:red_giros}). The translational subgroup of $Q$ 
is a cubic primitive lattice with vectors of length one.  $\nor(Q)$ is the group of symmetries of the set of triad rotation axes. It is a group of type $I\frac{4_{1}}{g}\overline{3}\frac{2}{d}$, and its
translational subgroup  is a body centered lattice generated by the vectors $(\pm \frac{1}{2}, \pm \frac{1}{2} ,\pm \frac{1}{2})$. 

\begin{figure}[htb]
\[
\begin{array}{ccccc}
&&
\includegraphics[width=3cm]{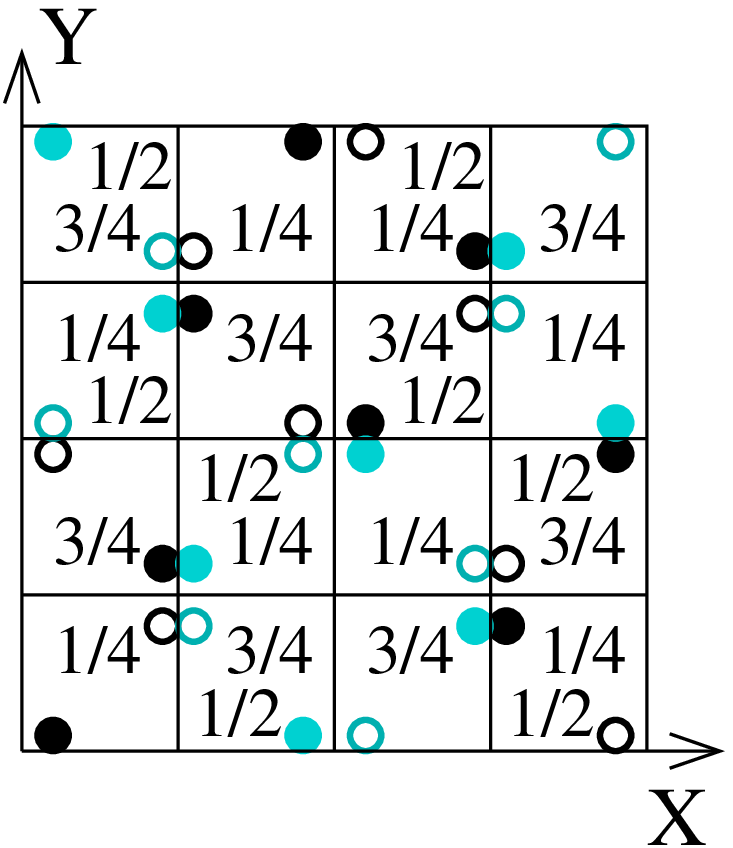}
&\\

&&
\nor(Q)=I\frac{4_{1}}{g}\overline{3}\frac{2}{d}
&\\
\\
&\swarrow&\downarrow&\searrow\\
\\
\includegraphics[width=3cm]{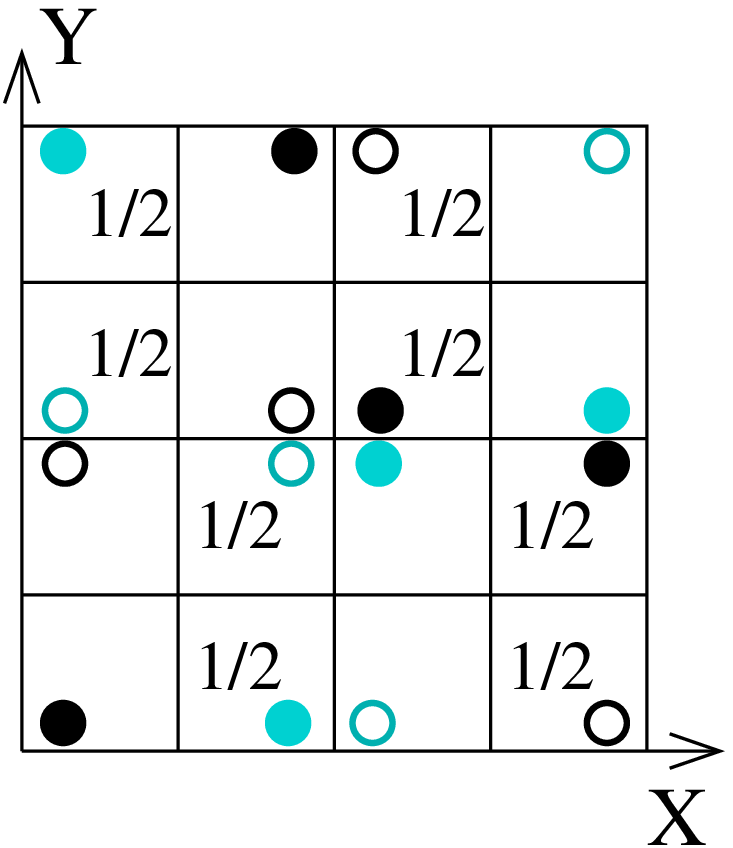}
&&
\includegraphics[width=3cm]{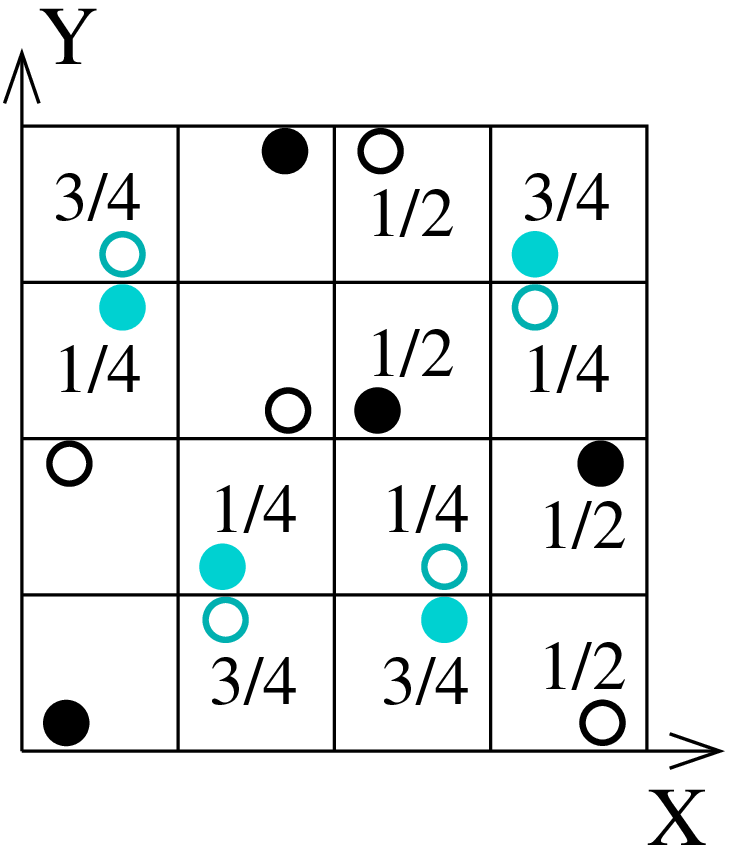}
&&
\includegraphics[width=3cm]{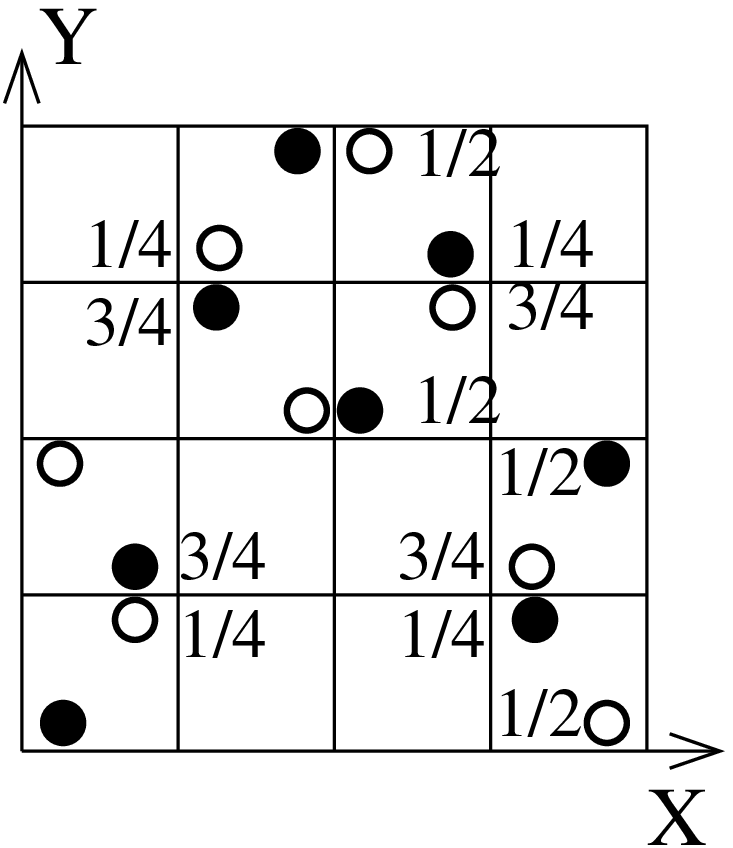}
\\

I\frac{2}{g}\overline{3}
&&
I\overline{4}3d
&&
I4_{1}32
\\

\\
\downarrow&\searrow&\downarrow&\swarrow&\downarrow\\
\\
\includegraphics[width=3cm]{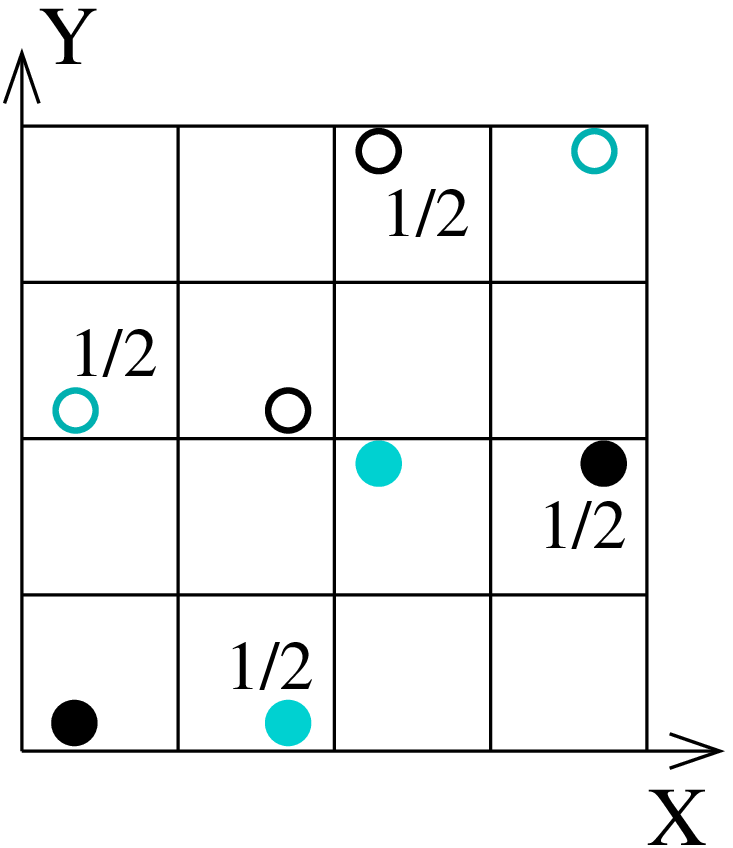}
&&
\includegraphics[width=3cm]{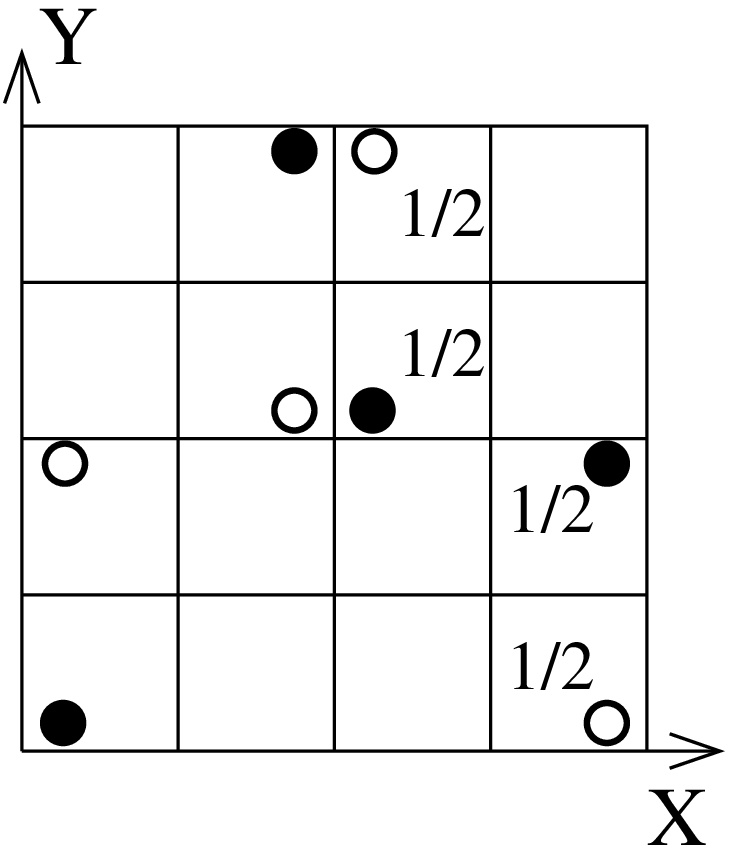}
&&
\includegraphics[width=3cm]{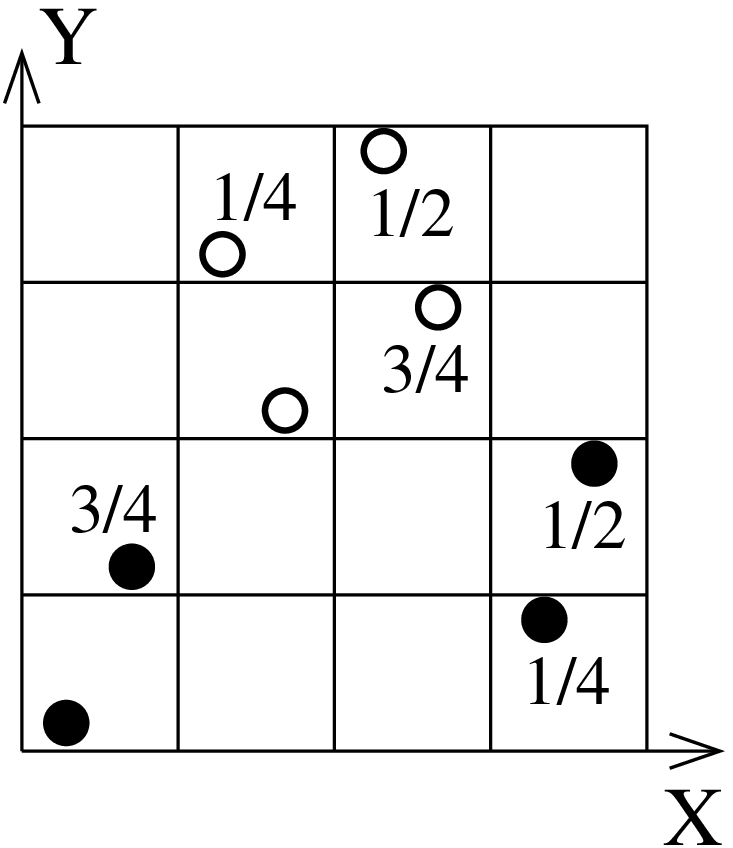}
\\

P\frac{2_{1}}{a}\overline{3}
&&
I2'3
&&
P4_{1}32
\\
\\
&\searrow&\downarrow&\swarrow&\\
\\
&&
\includegraphics[width=3cm]{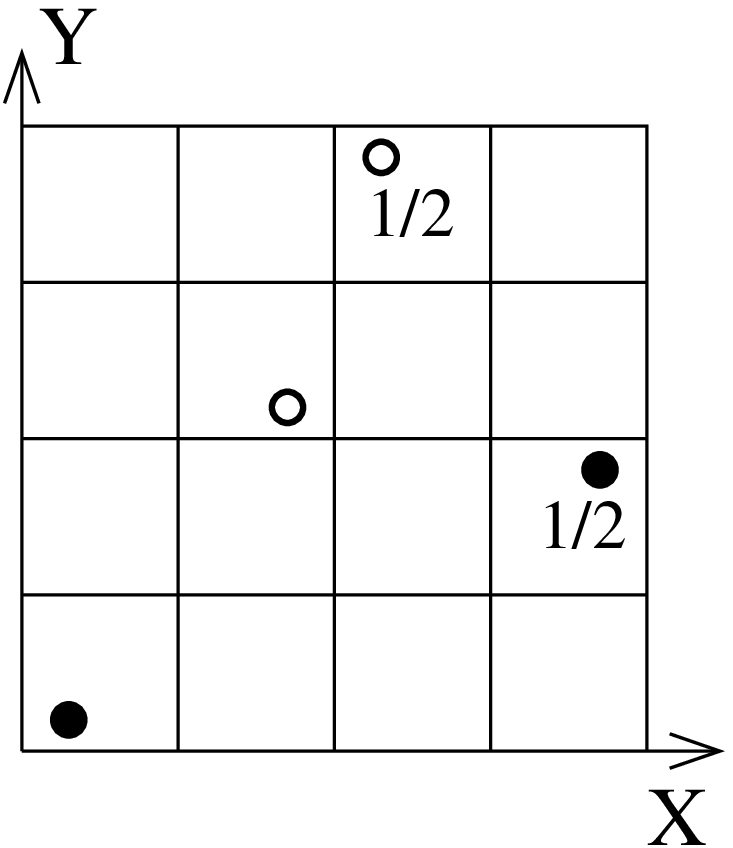}
&&
\\

&&
Q=P2_{1}3
&&
\\

\end{array}
\]
\caption{The eight quarter groups.
\label{fig:symmetries}}
\end{figure}

$\nor(Q)/Q$ has order 8. In fact, it is isomorphic to the dihedral group $D_{8}$. Hence there are 8 groups between $\nor(Q)$ and $Q$, including both. Their lattice  is drawn in Figure~\ref{fig:symmetries}.   
In Figure~\ref{fig:symmetries}, we represent them graphically with the conventions of  \cite{Lockwood-1978}: for each group,  the intersection of a generic orbit with the primitive cell $[0,1]^{3}$ is considered. Because all the groups contain the triad rotation on the diagonal axis $x=y=z$, only a third of the orbit is necessary in order to describe the group. E.~g., the orbit with respect to the subgroup that sends horizontal planes to horizontal planes. This is what is drawn in the figures,  projected over the $XY$ plane. 

The full dot near the origin represents a base point $(a,b,c)$ with $0 \le a,b,c \le 1/4$. Next to each of the other orbit points there is a number $h\in\{0,1/4,1/2,3/4\}$, omitted whenever it equals $0$. If an orbit point is represented as a ``full dot'', then its $z$-coordinate equals  $c+h$. If it is represented as an ``empty dot'', then it equals $h-c$ (or $1-c$, if $h$ is zero). 
For extra understanding of the group, each dot is drawn dark or light depending on whether it is obtained from the base point by an orientation preserving or an orientation reversing isometry in $G$.

\subsection{Outline of our method}
\label{sec:outline}

 Let $G$ be one of the eight quarter groups and let $p\in \R^{3}$ be a base point for an orbit $Gp$, so that the Dirichlet stereohedron we want to study is the closed Voronoi region $\Vor_{Gp}(p)$. There is no loss of generality in assuming that $p$ has trivial stabilizer, since Koch \cite{Koch73,Koch84} completely classified cubic orbits with non-trivial stabilizer (more generally, those with less than 3 degrees of freedom) and showed that they produce Dirichlet stereohedra with  at most 23 facets.
 
 If $p$ has trivial stabilizer, then a small perturbation of its coordinates can only increase the number of facets \cite[Lemma 3.1]{SabariegoSantos-full}, so we assume $p$ to be sufficiently generic.
 
Our method is as follows:
We consider a certain tessellation $\auxtes$ of the three-dimensional Euclidean space, which we call the {\em auxiliary tessellation}. For each tile $T$ of $\auxtes$ we call \emph{extended Voronoi region} of $T$ any region $\VorExt_{G}{T}$
with the following property:
\[
\forall q\in T, \qquad \Vor_{Gq}(q) \subseteq \VorExt_{G}(T).
\]
The extended Voronoi region is not uniquely defined. Part of our work is to compute one that is as small as possible, for each tile $T$.

Now, if $T_0$ is the tile of $\auxtes$ that contains the base point $p$, we call \emph{influence region of $T_0$} (or, of $p$) the union of all the tiles of $\auxtes$ whose extended Voronoi regions meet $\VorExt_{G}(T_0)$ in their interiors. 
We denote it by $\Infl_{G}(T_0)$. Observe that, strictly speaking, $\Infl_{G}(T_0)$ depends not only on $G$, but also on $\auxtes$ and on our particular choice of extended Voronoi regions.
It is easy to prove:

\begin{theorem}
\label{thm:influence}
If $T_0$ is the tile containing $p$, then all the
elements of $Gp$ that are 
neighbors of $p$ in the Voronoi diagram $\Vor_{Gp}(p)$ 
lie in $\Infl_{G}(T_0)$.
\end{theorem} 

\begin{proof}
See \cite{Bochis-1999}.
\end{proof}

That is, to get a bound on the number of facets of Dirichlet stereohedra for base points in $T_0$ we only need to count, or bound, $|Gp \cap \Infl_{G}(T_0)|$.

The above ``sketch of a method'' is impossible to implement directly. It involves the computation of an infinite number of extended Voronoi regions and influence regions. The key property that allows us to convert this into a finite computation is the following easy lemma:

 \begin{lemma}
\label{lem:normalizer}
Let $\nor(G)$ denote the normalizer of a crystallographic group $G$.
 If $T_1$ and $T_2=\rho T_1$ are tiles related by a transformation $\rho\in \nor(G)$, then 
 $\VorExt_{G}(T_2)=\rho \VorExt_{G}(T_{1})$ and
 $\Infl_{G}(T_2)=\rho \Infl_{G}(T_{1})$.
 \end{lemma}
 
 \begin{proof}
 It follows from the fact that $\Vor_{G \rho(q)}(\rho(q)) = \rho (\Vor_{Gq}(q))$ for every point $q$, if $\rho$ is in the normalizer $\nor(G)$ of $G$.
 \end{proof}

So, what we do is to use an auxiliary tessellation whose tiles lie in a small number of classes modulo the normalizer of $G$. Only this number of influence regions needs to be computed.

\subsection{Computation of extended Voronoi regions}
\label{sec:cut}
The basic idea for our computation of extended Voronoi regions is that each translation and rotation present in $G$ implies that a certain region of space can be excluded from $\VorExt_{G}(T_{0})$. 


\begin{lemma}
\label{lemma:VorExt-rotation}
Let $T_{0}$ be a convex domain and let $\rho$ be a rotation in $G$ with axis $\ell$ and rotation angle $\alpha$. 
Assume that $\ell$ does not intersect $T_0$. Let $H_1$ and $H_2$ be the two support half-planes of $T_0$ with border in  $\ell$. Let $H'_1$ and $H'_2$ be the half-planes obtained by rotating  $H_1$ and $H_2$ ``away from $T_0$'' with angles $\pm\alpha/2$ and axis $\ell$. 

Then, for every $p\in T_0$, the Dirichlet region $\Vor_{Gp}(p)$ is contained in the (perhaps non-convex) dihedral region bounded by $H'_1\cup H'_2$ and containing $T_0$.
\end{lemma}

\begin{figure}
\[
\begin{array}{cc}
      \includegraphics[width=5.8cm]{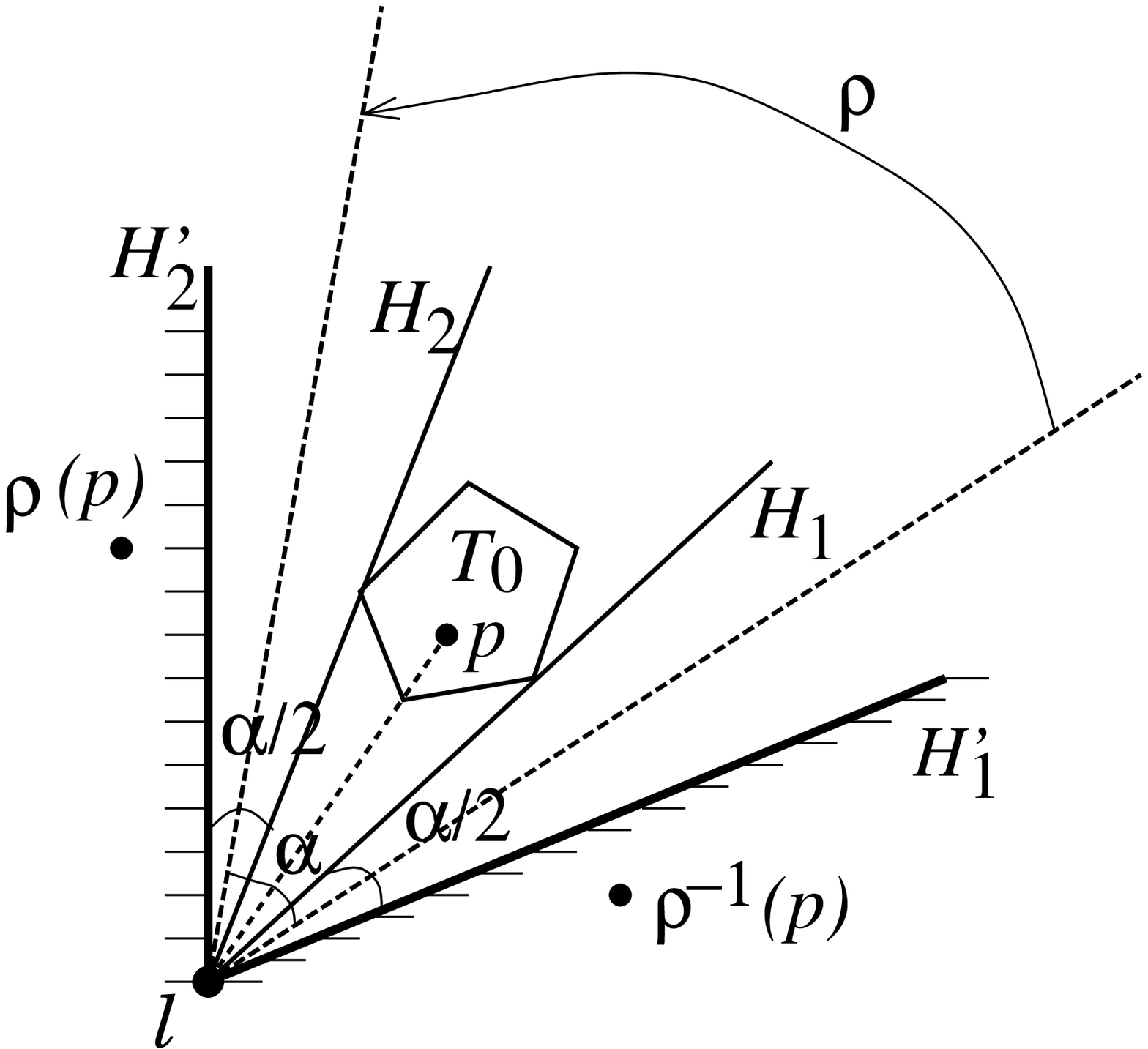}
    &
     \includegraphics[width=5.8cm]{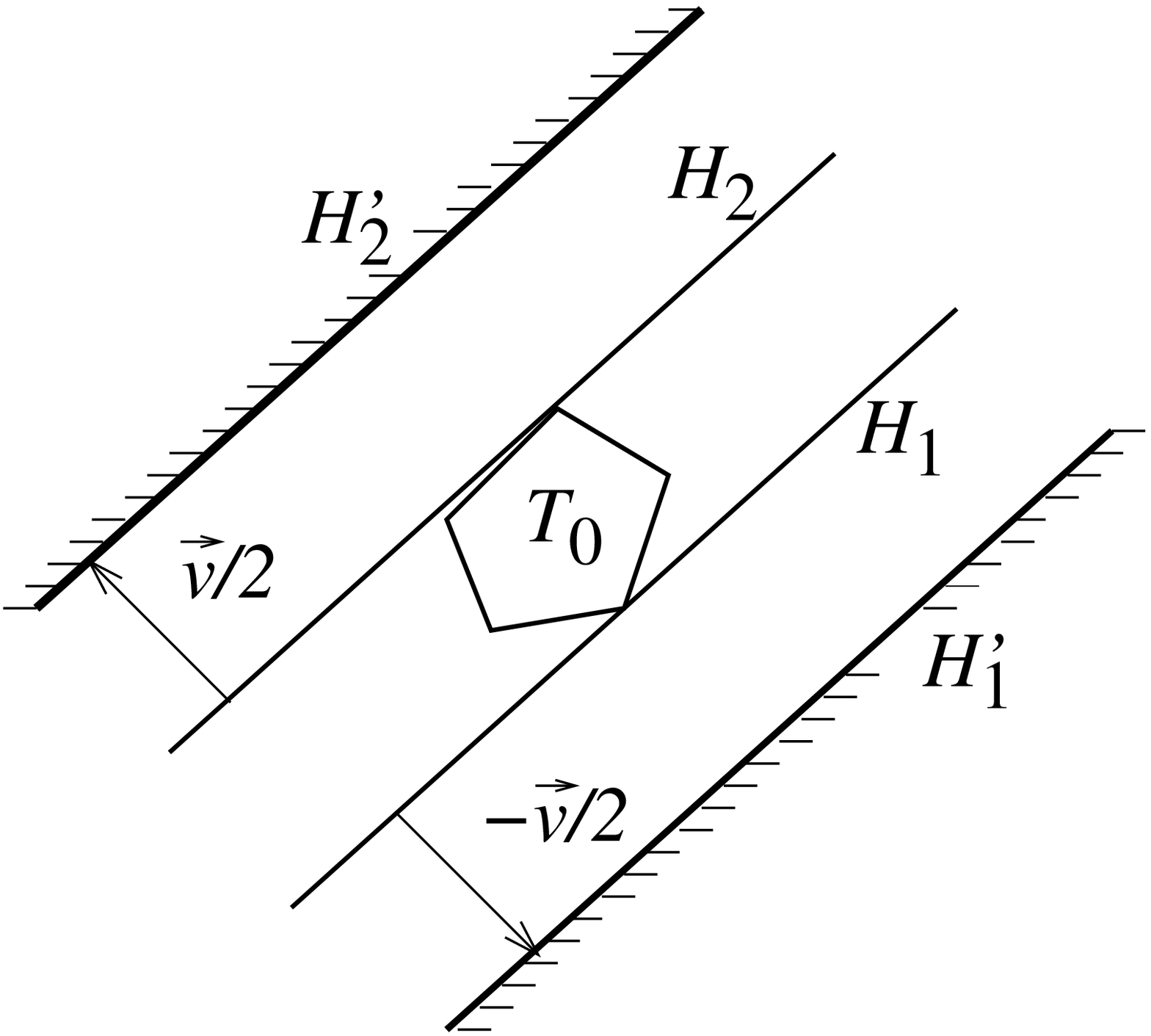}
    \end{array}
    \]
      \caption{Situations in the Lemmas~\ref{lemma:VorExt-rotation} and \ref{lemma:VorExt-translation}}
     \label{fig:VorExtRotTrans}
\end{figure}

\begin{proof}
In fact, something stronger is true: $\Vor_{Gp}(p)$ is contained in the dihedral sector with axis in $\ell$, angle $\alpha$, and centered at $p$. This follows from the fact that the two facets of this dihedron are parts of the bisectors of $p$ and $\rho p$ and of $p$ and $\rho^{-1} p$, respectively. See the left part of Figure~\ref{fig:VorExtRotTrans}.
\end{proof}

\begin{lemma}
\label{lemma:VorExt-translation}
Let $T_{0}$ be a convex domain and let $\vec{v}$ be the vector of some translation in $G$. 
Let $H_1$ and $H_2$ be the two support planes of $T_0$ orthogonal to $\vec{v}$.
Let $H'_1$ and $H'_2$ be the planes obtained translating  $H_1$ and $H_2$ ``away from $T_0$'' by the vectors $\pm\vec{v}/2$.

Then, for every $p\in T_0$, the Dirichlet region $\Vor_{Gp}(p)$ is contained in the 
strip between $H'_1$ and $H'_2$.
\end{lemma}

\begin{proof}
Similar to the previous one. See the right part of  Figure~\ref{fig:VorExtRotTrans}.
\end{proof}

\begin{figure}
  \begin{center}
      \includegraphics[width=5.5cm]{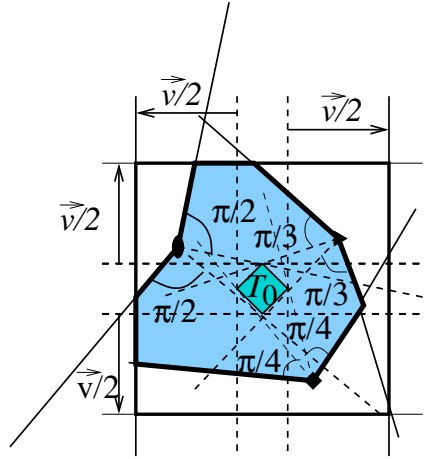}
    \end{center}
      \caption{Extended Voronoi region induced by three rotations and two translations}
     \label{fig:VorExt1}
 \end{figure}

 So, our method for constructing extended Voronoi regions consists in identifying a certain number of rotations and translations in $G$ and taking as $\VorExt_{G}(T_0)$ the intersection of the regions allowed in Lemmas~\ref{lemma:VorExt-rotation} and~\ref{lemma:VorExt-translation}.
 Figure~\ref{fig:VorExt1}  shows an example of what we mean, in which we consider three rotations and two translations in the plane.
 
The problem with this approach is that the regions obtained in Lemma~\ref{lemma:VorExt-rotation} may not be convex. Observe that we are going to use more than 250 rotations (see Tables~\ref{table:axes-triad-rotations}, \ref{table:coordinate-diad-rotations}, \ref{table:diagonal-diad-rotations} and~\ref{table:diagonal-P4132}). So, in principle, to compute an extended Voronoi region we would need to intersect 250 non-convex regions, which is an extremely hard computational problem (see, e.~g.,~\cite{Erickson-polyhedra}).

To avoid this we do the following, at the expense of getting a slightly larger (hence, for our purposes worse) extended Voronoi region: instead of obtaining the extended Voronoi region directly as an intersection, we obtain it as a union of tiles of the same auxiliary tessellation $\auxtes$. That is to say: 
\begin{enumerate}
\item We first identify a finite (but large) population of tiles in $\auxtes$ that form themselves an extended Voronoi region (that is, which are guaranteed to contain $\Vor_{Gp}(p)$ for every $p\in T_0$).

\item We then process one by one the rotations and translations we are interested in, and at each step discard from our initial population the tiles that do not intersect the corresponding region of Lemma~\ref{lemma:VorExt-rotation} or Lemma~\ref{lemma:VorExt-translation}.
\end{enumerate}

\begin{figure}
    \begin{center}
      \includegraphics[width=8cm]{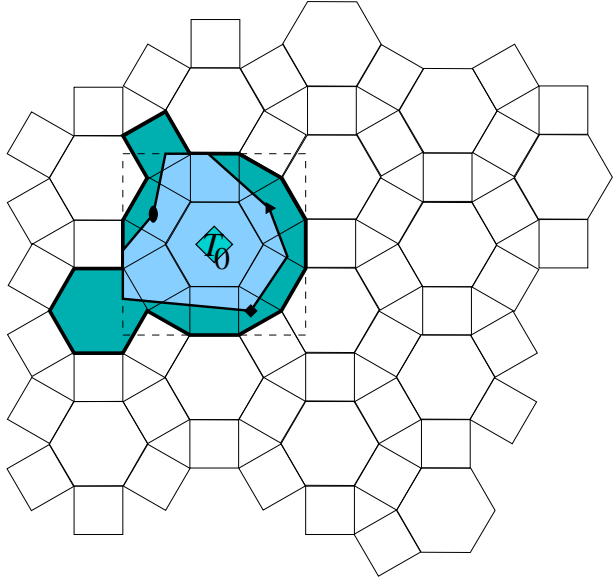}
    \end{center}
      \caption{Extended Voronoi region over a tessellation}
     \label{fig:VorExt2}
  \end{figure}

What we eventually obtain with this method is the set of tiles of $\auxtes$ that intersect the region of Figure~\ref{fig:VorExt1}. See Figure~\ref{fig:VorExt2}.

  The reader may wonder why step (2) above is computationally not so hard, since it again involves the intersection of non-convex regions. The reason is that now we intersect only two such regions at a time (one tile and one strip or dihedron) and do not need to intersect the resulting set with anything else. Also (although less important), we do not really need to compute an intersection but only to check its emptyness.

\section{Detailed discussion of the method}

\subsection{The auxiliary tessellation}

The first natural choice of an auxiliary tessellation would be a tessellation by fundamental domains of the normalizer $\nor(G)$, or of any subgroup of it. In this way  all the tiles are equivalent: By Lemma~\ref{lem:normalizer}, in such a tessellation we would need to compute only one extended Voronoi region. 

\begin{figure}
    \begin{center}
      \includegraphics[width=5.5cm]{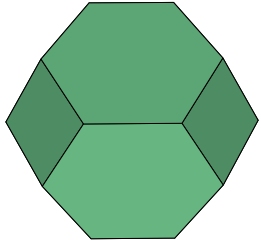}
    \end{center}
      \caption{A truncated octahedron, the Voronoi cell of the body centered cubic lattice}
     \label{fig:octahedron}
\end{figure}

Since we have $\nor(Q)\le \nor(G)$ for every quarter group, as a first step let us compute a fundamental domain of $\nor(Q)$ (later, we will subdivide it further). In order to get a fundamental domain with some symmetries, we start with the Voronoi region of one of its degenerate orbits, namely the orbit with base point at the origin $(0,0,0)$.
This orbit is a body centered cubic lattice of side $1/2$. Its Voronoi cell is a truncated octahedron (see Figure~\ref{fig:octahedron}) whose 14 facets have supporting inequalities
\[
\pm x \le  \frac{1}{4},\quad
\pm y \le  \frac{1}{4},\quad
\pm z \le  \frac{1}{4},\quad
\pm x \pm y \pm z \le \frac{3}{8}.
\]

The  24 vertices of the Voronoi cell are 
\[
\begin{matrix}
\left(\pm \frac{1}{4}, \pm \frac{1}{8}, 0\right),&
\left(\pm \frac{1}{4}, 0, \pm \frac{1}{8}\right),&
\left(0, \pm \frac{1}{4}, \pm \frac{1}{8}\right),\\
\\
\left(\pm \frac{1}{8}, \pm \frac{1}{4}, 0\right),&
\left(\pm \frac{1}{8}, 0, \pm \frac{1}{4}\right),&
\left(0, \pm \frac{1}{8}, \pm \frac{1}{4}\right).
\end{matrix}
\]

The stabilizer of this truncated octahedron in $\nor(Q)$ has order six. It is generated by the central symmetry (inversion) around the origin and the rotation of order three with axis  in the line $x=z=y$. Hence, one way of obtaining a fundamental domain for $\nor(Q)$ is to intersect the truncated octahedron with  
a sector of  angle $\pi/3$ with edge in that line. There are several inequivalent ways of doing it, but we choose the sector defined by the inequalities $x\le y\le z$.

\begin{figure}
\[
\begin{array}{c|c}
      \includegraphics[width=5.9cm]{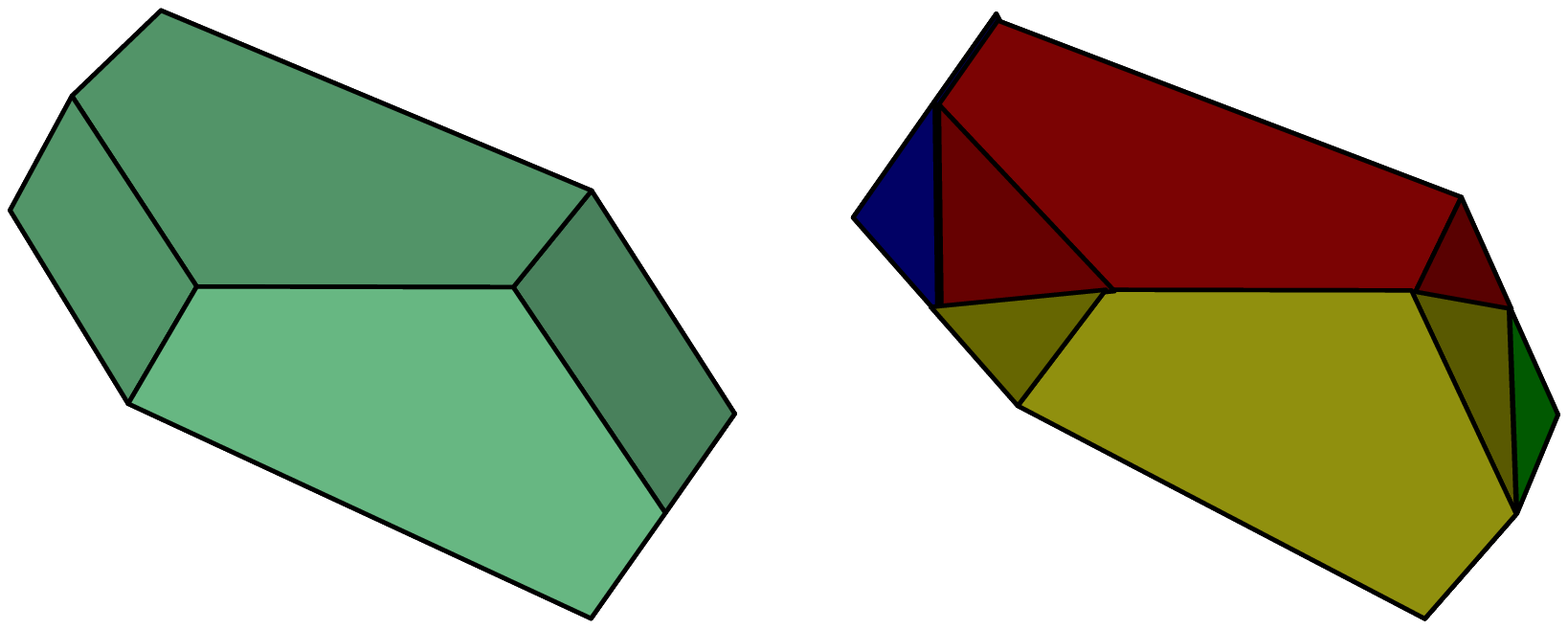}
    &
     \includegraphics[width=5.9cm]{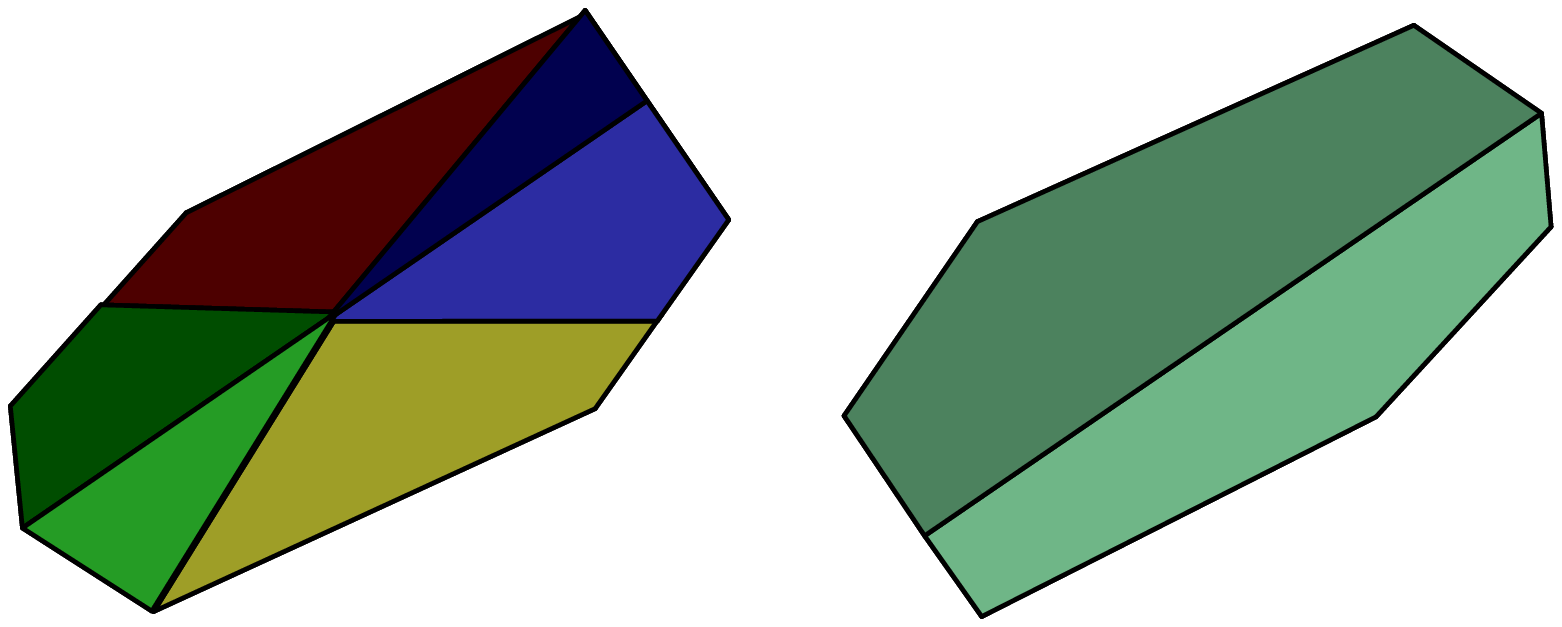}
    \end{array}
    \]
      \caption{Two views of the fundamental domain of $\nor(Q)$ and its subdivision into four prototiles, $A_{0}$, $B_{0}$, $C_{0}$ and $D_{0}$}
     \label{fig:domains}
     \end{figure}

In this way, we obtain as fundamental domain of $\nor(Q)$ the polytope with the following 12 vertices:
\[
\begin{matrix}
\left(   \frac{1}{8},   \frac{1}{8},   \frac{1}{8}         \right) \\
\\
\left(   \frac{1}{4},   \frac{1}{16},   \frac{1}{16}         \right) \qquad
\left(   \frac{1}{4},   \frac{1}{8},   0         \right) \qquad
\left(   \frac{3}{16},   \frac{3}{16},   0         \right) \\
\\
\left(   \frac{1}{4},   \frac{-1}{16},   \frac{-1}{16}         \right) \qquad
\left(   \frac{1}{4},   0,   -\frac{1}{8}         \right) \qquad
\left(   \frac{1}{8},   0,   \frac{-1}{4}         \right) \qquad
\left(   \frac{1}{16},   \frac{1}{16},   \frac{-1}{4}         \right) \\
\\
\left(   0, \frac{-3}{16},   \frac{-3}{16}          \right) \qquad
\left(   0, \frac{-1}{8},    \frac{-1}{4}           \right) \qquad
\left(   \frac{-1}{16},   \frac{-1}{16},   \frac{-1}{4}         \right) \\
\\
\left(   \frac{-1}{8},   \frac{-1}{8},   \frac{-1}{8}         \right) \\
\end{matrix}
\]
This polytope is depicted on the left side of Figure~\ref{fig:domains}. It has eight facets: the two ``big ones'' in the planes $x=y$ and $y=z$ are pentagons and the other six are parts of facets of the truncated octahedron: two halves-of-squares, two halves-of-hexagons and two sixths-of-hexagons.

As we said above, to get better bounds, we further subdivide this fundamental domain of $\nor(Q)$.
We do it cutting with the coordinate planes $z=0$, $y=0$ and $x=0$. We get the following four sub-polytopes  $A_{0}$, $B_{0}$, $C_{0}$ and $D_{0}$ (see Figure~\ref{fig:domains}):
\[
A_{0}:= \conv\hspace{-2.5pt}\left(\hspace{-2.5pt}
(0,0,0), 
\left(   \frac{1}{8},   \frac{1}{8},   \frac{1}{8}         \right),
\left(   \frac{1}{4},   \frac{1}{16},   \frac{1}{16}         \right),
\left(   \frac{1}{4},   \frac{1}{8},   0         \right),
\left(   \frac{3}{16},   \frac{3}{16},   0         \right),
\left(   \frac{1}{4},   0,   0         \right) 
\right)
\]

\begin{eqnarray*}
B_{0} &:=& \conv\left(
\left(0,0,0\right),
\left(   \frac{1}{4},   0, 0         \right) ,
\left(   \frac{1}{4},   \frac{1}{8},   0         \right) ,
\left(   \frac{3}{16},   \frac{3}{16},   0         \right) , \right. \\
& &
\left.
\left(   \frac{1}{4},   0,   \frac{-1}{8}         \right) ,
\left(   \frac{1}{8},   0,   \frac{-1}{4}         \right) ,
\left(   \frac{1}{16},   \frac{1}{16},   \frac{-1}{4}         \right) ,
\left( 0, 0, \frac{-1}{4}    \right)
\right)
\end{eqnarray*}

\begin{eqnarray*}
C_{0} &:= &\conv\left(
\left(0,0,0\right),
\left(   \frac{1}{4},   0, 0         \right) ,
\left(   \frac{1}{4},   \frac{-1}{16},   \frac{-1}{16}         \right) ,
\left(   \frac{1}{4},   0,   \frac{-1}{8}         \right) , \right. \\
& &
\left.
 \left(   \frac{1}{8},   0,   \frac{-1}{4}         \right) ,
\left(   0, \frac{-3}{16},   \frac{-3}{16}          \right) ,
\left(   0, \frac{-1}{8},    \frac{-1}{4}           \right) ,
\left( 0, 0, \frac{-1}{4}   \right)
\right)
\end{eqnarray*}

\begin{eqnarray*}
D_{0} &:= &\conv\left(
\left(0,0,0\right),
\left(   0, \frac{-3}{16},   \frac{-3}{16}          \right) ,
\left(   \frac{-1}{16},   \frac{-1}{16},   \frac{-1}{4}         \right) , \right. \\
& &
\left.
\left(   0, \frac{-1}{8},    \frac{-1}{4}           \right) ,
\left(   \frac{-1}{8},   \frac{-1}{8},   \frac{-1}{8}         \right),
\left( 0, 0, \frac{-1}{4}     \right) 
\right)
\end{eqnarray*}

\begin{definition}
\label{defi:auxtes}
From now on, the tessellation $\auxtes$ of $\R^3$ obtained letting $\nor(Q)$ act on the four  ``prototiles'' $A_0$, $B_0$, $C_0$, and $D_0$ is called \emph{auxiliary tessellation}. 
\end{definition}

\subsection{The initial population of tiles}

As said in Section~\ref{sec:cut}, we compute extended Voronoi regions starting with a finite number of tiles of $\auxtes$ whose union in guaranteed to contain $\Vor_{Gp}(p)$ for every $p\in T_{0}$, and then eliminating as many of them as we can. As initial population of  tiles, we take all the ones contained in the cube $[-1,1]^3$. That this is enough follows from the fact that the four prototiles are contained in the cube $[-1/4,1/4]^3$ and our group $Q$ contains the translations of length $1$ in the three coordinate directions. Lemma~\ref{lemma:VorExt-translation} then implies that the extended Voronoi region (of any of the prototiles) is contained in the cube $[-3/4,3/4]^3$.

Since the four prototiles together form a fundamental domain of $\nor(Q)$ and since the density of (any generic orbit of) $\nor(Q)$ is 96 points per unit cube, in the initial population we have $96 \times 8 = 768$ tiles of each of the four types, and  $3\, 072$ tiles in total.

In our implementation, each tile of $\auxtes$ is characterized by its ``type'' ($A$, $B$, $C$ or $D$) and the transformation $\rho\in \nor(Q)$ that sends the corresponding prototile ($A_0$, $B_0$, $C_0$ or $D_0$) to it. We compute the initial tiles of type, say, $A$, in the initial tessellation as follows:

\begin{itemize}

\item From the graphic representation of $\nor(Q)$ in Figure~\ref{fig:symmetries} we read the 32 isometries in $\nor(Q)$ that  send horizontal planes to horizontal planes
and send a base point $p\in[0,1/4]^3$ to lie in the unit cube $[0,1]^3$.
The matrices of these isometries are shown in Table~\ref{table:motionsA}. For future reference, these 32
matrices appear in eight subsets $\mathcal{C}_{i}$, $i=1,\dots, 8$, representing the cosets modulo $Q$. Table~\ref{table:motions-groups} shows which cosets of isometries form each of the eight quarter groups. 


\begin{table}
\footnotesize
\[
\begin{array}{|c|c|c|c|c|}
\hline
& \multicolumn{4}{c|}{\text{\normalsize Type {\em A}}} \\  
\hline
\text{\normalsize $\mathcal{C}_{1}$} &
\begin{matrix}
1 & 0 & 0 & 0 \\
0 & 1 & 0 & 0 \\
0 & 0 & 1 & 0 \\
0 & 0 & 0 & 1 
 \end{matrix} 
%
 & 
  
\begin{matrix}
1 & 0 & 0 & 0 \\
1 & -1 & 0 & 0 \\
1/2 & 0 & -1 & 0 \\
1/2 & 0 & 0 & 1 
\end{matrix}
 
 & 
  
\begin{matrix}
1 & 0 & 0 & 0 \\
1/2 & -1 & 0 & 0 \\
1/2 & 0 & 1 & 0 \\
1 & 0 & 0 & -1 
\end{matrix}
 
&
  
\begin{matrix}
1 & 0 & 0 & 0 \\
1/2 & 1 & 0 & 0 \\
1 & 0 & -1 & 0 \\
1/2 & 0 & 0 & -1 
\end{matrix}
 
 \\
 \hline
 \text{\normalsize $\mathcal{C}_{2}$}  &
 
\begin{matrix}
1 & 0 & 0 & 0 \\
3/4 & 0 & 1 & 0 \\
1/4 & -1 & 0 & 0 \\
1/4 & 0 & 0 & 1 
\end{matrix}
 
&
   
\begin{matrix}
1 & 0 & 0 & 0 \\
1/4 & 0 & -1 & 0 \\
1/4 & 1 & 0 & 0 \\
3/4 & 0 & 0 & 1 
\end{matrix}
 
 & 
  
\begin{matrix}
1 & 0 & 0 & 0 \\
3/4 & 0 & -1 & 0 \\
3/4 & -1 & 0 & 0 \\
3/4 & 0 & 0 & -1 
\end{matrix}
 
& 
  
\begin{matrix}
1 & 0 & 0 & 0 \\
1/4 & 0 & 1 & 0 \\
3/4 & 1 & 0 & 0 \\
1/4 & 0 & 0 & -1 
\end{matrix}
 
\\
\hline
\text{\normalsize $\mathcal{C}_{3} $} &
 
\begin{matrix}
1 & 0 & 0 & 0 \\
1 & -1 & 0 & 0 \\
0 & 0 & 1 & 0 \\
1/2 & 0 & 0 & -1 
\end{matrix}
 
&
  
\begin{matrix}
1 & 0 & 0 & 0 \\
0 & 1 & 0 & 0 \\
1/2 & 0 & -1 & 0 \\
1 & 0 & 0 & -1 
\end{matrix}
 
 & 
 
\begin{matrix}
1 & 0 & 0 & 0 \\
1/2 & 1 & 0 & 0 \\
1/2 & 0 & 1 & 0 \\
1/2 & 0 & 0 & 1 
\end{matrix}
 
&
 
\begin{matrix}
1 & 0 & 0 & 0 \\
1/2 & -1 & 0 & 0 \\
1 & 0 & -1 & 0 \\
0 & 0 & 0 & 1 
\end{matrix}
 
\\
 \hline
 \text{\normalsize $\mathcal{C}_{4}$}  &
  
\begin{matrix}
1 & 0 & 0 & 0 \\
1/4 & 0 & -1 & 0 \\
1/4 & -1 & 0 & 0 \\
1/4 & 0 & 0 & -1 
\end{matrix}
 
 & 
 
\begin{matrix}
1 & 0 & 0 & 0 \\
3/4 & 0 & 1 & 0 \\
1/4 & 1 & 0 & 0 \\
3/4 & 0 & 0 & -1 
\end{matrix}
 
 & 
  
\begin{matrix}
1 & 0 & 0 & 0 \\
1/4 & 0 & 1 & 0 \\
3/4 & -1 & 0 & 0 \\
3/4 & 0 & 0 & 1 
\end{matrix}
 
&
  
\begin{matrix}
1 & 0 & 0 & 0 \\
3/4 & 0 & -1 & 0 \\
3/4 & 1 & 0 & 0 \\
1/4 & 0 & 0 & 1 
\end{matrix}
 
 \\
 \hline
\text{\normalsize $\mathcal{C}_{5} $} &
 
\begin{matrix}
1 & 0 & 0 & 0 \\
1/2 & 1 & 0 & 0 \\
0 & 0 & 1 & 0 \\
1 & 0 & 0 & -1 
\end{matrix}
 
& 
  
\begin{matrix}
1 & 0 & 0 & 0 \\
1/2 & -1 & 0 & 0 \\
1/2 & 0 & -1 & 0 \\
1/2 & 0 & 0 & -1 
\end{matrix}
 
 &
 
\begin{matrix}
1 & 0 & 0 & 0 \\
1 & -1 & 0 & 0 \\
1/2 & 0 & 1 & 0 \\
0 & 0 & 0 & 1 
\end{matrix}
 
 &
   
\begin{matrix}
1 & 0 & 0 & 0 \\
0 & 1 & 0 & 0 \\
1 & 0 & -1 & 0 \\
1/2 & 0 & 0 & 1 
\end{matrix}
 
\\
 \hline
 \text{\normalsize $\mathcal{C}_{6}$}  &
 
\begin{matrix}
1 & 0 & 0 & 0 \\
1/4 & 0 & 1 & 0 \\
1/4 & -1 & 0 & 0 \\
3/4 & 0 & 0 & -1 
\end{matrix}
 
 &
  
\begin{matrix}
1 & 0 & 0 & 0 \\
3/4 & 0 & -1 & 0 \\
1/4 & 1 & 0 & 0 \\
1/4 & 0 & 0 & -1 
\end{matrix}
 
& 
 
\begin{matrix}
1 & 0 & 0 & 0 \\
1/4 & 0 & -1 & 0 \\
3/4 & -1 & 0 & 0 \\
1/4 & 0 & 0 & 1 
\end{matrix}
 
 & 
 
\begin{matrix}
1 & 0 & 0 & 0 \\
3/4 & 0 & 1 & 0 \\
3/4 & 1 & 0 & 0 \\
3/4 & 0 & 0 & 1 
\end{matrix}
 
\\
 \hline
 \text{\normalsize $\mathcal{C}_{7}$}  &
 
\begin{matrix}
1 & 0 & 0 & 0 \\
1/2 & -1 & 0 & 0 \\
0 & 0 & 1 & 0 \\
1/2 & 0 & 0 & 1 
\end{matrix}
 
 &
  
\begin{matrix}
1 & 0 & 0 & 0 \\
1/2 & 1 & 0 & 0 \\
1/2 & 0 & -1 & 0 \\
0 & 0 & 0 & 1 
\end{matrix}
 
& 
  
\begin{matrix}
1 & 0 & 0 & 0 \\
0 & 1 & 0 & 0 \\
1/2 & 0 & 1 & 0 \\
1/2 & 0 & 0 & -1 
\end{matrix}
 
 & 
 
\begin{matrix}
1 & 0 & 0 & 0 \\
1 & -1 & 0 & 0 \\
1 & 0 & -1 & 0 \\
1 & 0 & 0 & -1 
\end{matrix}
 
 \\
 \hline
\text{\normalsize $\mathcal{C}_{8} $} &
 
\begin{matrix}
1 & 0 & 0 & 0 \\
3/4 & 0 & -1 & 0 \\
1/4 & -1 & 0 & 0 \\
3/4 & 0 & 0 & 1 
\end{matrix}
 
 & 
  
\begin{matrix}
1 & 0 & 0 & 0 \\
1/4 & 0 & 1 & 0 \\
1/4 & 1 & 0 & 0 \\
1/4 & 0 & 0 & 1 
\end{matrix}
 
 & 
 
\begin{matrix}
1 & 0 & 0 & 0 \\
3/4 & 0 & 1 & 0 \\
3/4 & -1 & 0 & 0 \\
1/4 & 0 & 0 & -1 
\end{matrix}
 
&
 
\begin{matrix}
1 & 0 & 0 & 0 \\
1/4 & 0 & -1 & 0 \\
3/4 & 1 & 0 & 0 \\
3/4 & 0 & 0 & -1 
\end{matrix}
 
 \\
\hline
 \end{array}
 \]
 \caption{The $32$ matrices of the isometries which send horizontal planes to horizontal planes for $A_{0}$}
 \label{table:motionsA}
 \end{table}


\begin{table}
\[
\begin{array}{|c|c|}
 \hline
\text{Quarter groups} & \text{Cosets}\\
\hline
I\frac{4_{1}}{g}\overline{3}\frac{2}{d} & \mathcal{C}_{1},\dots, \mathcal{C}_{8} \\
\hline
I4_{1}32 & \mathcal{C}_{1}, \mathcal{C}_{2}, \mathcal{C}_{3}, \mathcal{C}_{4}  \\ 
\hline
I\overline{4}3d & \mathcal{C}_{1}, \mathcal{C}_{3}, \mathcal{C}_{6}, \mathcal{C}_{8}  \\
\hline
I\frac{2}{g}\overline{3} &\mathcal{C}_{1}, \mathcal{C}_{3}, \mathcal{C}_{5}, \mathcal{C}_{7}  \\
\hline
P4_{1}32 &\mathcal{C}_{1}, \mathcal{C}_{2}  \\
\hline
I2'3 & \mathcal{C}_{1}, \mathcal{C}_{3}\\
\hline
P\frac{2_{1}}{a}\overline{3} & \mathcal{C}_{1}, \mathcal{C}_{7} \\
\hline
P2_{1}3 & \mathcal{C}_{1} \\
\hline
\end{array}
\]
\caption{Cosets of isometries in each quarter group}
\label{table:motions-groups}
\end{table} 

\item To these 32 isometries we apply the triad rotation with axis $x=y=z$, to obtain the 96 elements of $\nor(Q)$ that send a point $p\in A_0\subset[0,1/4]^3$ to lie in $[0,1]^3$.

\item Finally, we compose these 96 with the translations by the vectors $(0,0,0)$, $(-1,0,0)$, $(0,-1,0)$, $(0,0,-1)$, $(-1,-1,0)$, $(-1,0,-1)$, $(0,-1,-1)$ and $(-1,-1,-1)$ in order to get all the tiles of type $A$ contained in
$[-1,1]^3$.
\end{itemize}

With the tiles of the other types, $B$, $C$ and $D$, we work the same way, except for the fact that the prototiles are now not in $[0,1]^3$. For example: 
$B_{0}$ is in the unit cube $[0,1]\times[0,1]\times[-1,0]$. Still, to build our initial population of $B$-tiles, 
we want to start with the 32 horizontal transformations that send $B_0$ to lie in $[0,1]^3$. These are the same 32 of Table~\ref{table:motionsA} except some of them need to be composed with a suitable integer translation. The matrices that need modification are shown in Tables~\ref{table:motionsB},~\ref{table:motionsC},  and~\ref{table:motionsD}, for types $B$, $C$ and $D$, respectively.
As in case $A$, once we have got these isometries, we apply to them the triad rotations with axis $x=y=z$ and the translations  by the vectors $(0,0,0)$, $(-1,0,0)$, $(0,-1,0)$, $(0,0,-1)$, $(-1,-1,0)$, $(-1,0,-1)$, $(0,-1,-1)$ and $(-1,-1,-1)$ to get an initial population of $96\times 8 = 768$ tiles of each type.


\begin{table}
\small
\[
\begin{array}{|c|c|c|c|c|}
\hline
& \multicolumn{4}{c|}{\text{\normalsize  Type {\em B}}} \\  
\hline
\text{\normalsize $\mathcal{C}_{1}$} &
\begin{matrix}
1 & 0 & 0 & 0 \\
0 & 1 & 0 & 0 \\
0 & 0 & 1 & 0 \\
1 & 0 & 0 & 1 
\end{matrix}
 & 
 & 
 \begin{matrix}
1 & 0 & 0 & 0 \\
1/2 & -1 & 0 & 0 \\
1/2 & 0 & 1 & 0 \\
0 & 0 & 0 & -1 
\end{matrix}
&
 \\
 \hline
 %
 %
\text{\normalsize $\mathcal{C}_{3}$}  &
&
 \begin{matrix}
1 & 0 & 0 & 0 \\
0 & 1 & 0 & 0 \\
1/2 & 0 & -1 & 0 \\
0 & 0 & 0 & -1 
\end{matrix}
 & 
&
\begin{matrix}
1 & 0 & 0 & 0 \\
1/2 & -1 & 0 & 0 \\
1 & 0 & -1 & 0 \\
1 & 0 & 0 & 1 
\end{matrix}
\\
 \hline
 %
 %
\text{\normalsize $\mathcal{C}_{5}$}  &
\begin{matrix}
1 & 0 & 0 & 0 \\
1/2 & 1 & 0 & 0 \\
0 & 0 & 1 & 0 \\
0 & 0 & 0 & -1 
\end{matrix}
& 
 &
\begin{matrix}
1 & 0 & 0 & 0 \\
1 & -1 & 0 & 0 \\
1/2 & 0 & 1 & 0 \\
1 & 0 & 0 & 1 
\end{matrix}
 &
\\
 \hline
 %
 %
 \text{\normalsize $\mathcal{C}_{7}$}  &
 &
 \begin{matrix}
1 & 0 & 0 & 0 \\
1/2 & 1 & 0 & 0 \\
1/2 & 0 & -1 & 0 \\
1 & 0 & 0 & 1 
\end{matrix}
& 
 & 
\begin{matrix}
1 & 0 & 0 & 0 \\
1 & -1 & 0 & 0 \\
1 & 0 & -1 & 0 \\
0 & 0 & 0 & -1 
\end{matrix}
 \\
 \hline
 %
 %
 \end{array}
 \]
 \caption{The new matrices of the isometries which send horizontal planes to horizontal planes for $B_{0}$}
 \label{table:motionsB}
 \end{table}


\begin{table}
\small
\[
\begin{array}{|c|c|c|c|c|}
\hline
& \multicolumn{4}{c|}{\text{\normalsize  Type {\em C}}} \\  
\hline
\text{\normalsize $\mathcal{C}_{1}$} &
 
\begin{matrix}
1 & 0 & 0 & 0 \\
0 & 1 & 0 & 0 \\
1 & 0 & 1 & 0 \\
1 & 0 & 0 & 1 
\end{matrix}
 
 & 
%
%
 & 
  
\begin{matrix}
1 & 0 & 0 & 0 \\
1/2 & -1 & 0 & 0 \\
1/2 & 0 & 1 & 0 \\
0 & 0 & 0 & -1 
\end{matrix}
 
&
  
\begin{matrix}
1 & 0 & 0 & 0 \\
1/2 & 1 & 0 & 0 \\
0 & 0 & -1 & 0 \\
1/2 & 0 & 0 & -1 
\end{matrix}
 
 \\
 \hline
 %
 %
%
%
%
%
%
%
%
%
\text{\normalsize $\mathcal{C}_{3}$}  &
 
\begin{matrix}
1 & 0 & 0 & 0 \\
1 & -1 & 0 & 0 \\
1 & 0 & 1 & 0 \\
1/2 & 0 & 0 & -1 
\end{matrix}
 
&
  
\begin{matrix}
1 & 0 & 0 & 0 \\
0 & 1 & 0 & 0 \\
1/2 & 0 & -1 & 0 \\
0 & 0 & 0 & -1 
\end{matrix}
 
 & 
%
%
&
 
\begin{matrix}
1 & 0 & 0 & 0 \\
1/2 & -1 & 0 & 0 \\
0 & 0 & -1 & 0 \\
1 & 0 & 0 & 1 
\end{matrix}
 
\\
 \hline
 %
 %
%
%
%
%
%
%
%
%
\text{\normalsize $\mathcal{C}_{5}$}  &
 
\begin{matrix}
1 & 0 & 0 & 0 \\
1/2 & 1 & 0 & 0 \\
1 & 0 & 1 & 0 \\
0 & 0 & 0 & -1 
\end{matrix}
 
& 
%
%
 &
 
\begin{matrix}
1 & 0 & 0 & 0 \\
1 & -1 & 0 & 0 \\
1/2 & 0 & 1 & 0 \\
1 & 0 & 0 & 1 
\end{matrix}
 
 &
   
\begin{matrix}
1 & 0 & 0 & 0 \\
0 & 1 & 0 & 0 \\
0 & 0 & -1 & 0 \\
1/2 & 0 & 0 & 1 
\end{matrix}
 
\\
 \hline
 %
 %
%
%
%
%
%
%
%
%
 \text{\normalsize $\mathcal{C}_{7}$}  &
 
\begin{matrix}
1 & 0 & 0 & 0 \\
1/2 & -1 & 0 & 0 \\
1 & 0 & 1 & 0 \\
1/2 & 0 & 0 & 1 
\end{matrix}
 
 &
  
\begin{matrix}
1 & 0 & 0 & 0 \\
1/2 & 1 & 0 & 0 \\
1/2 & 0 & -1 & 0 \\
1 & 0 & 0 & 1 
\end{matrix}
 
& 
%
%
 & 
 
\begin{matrix}
1 & 0 & 0 & 0 \\
1 & -1 & 0 & 0 \\
0 & 0 & -1 & 0 \\
0 & 0 & 0 & -1 
\end{matrix}
 
 \\
 \hline
 %
 %
%
%
%
%
%
%
%
%
\end{array}
 \]
 \caption{The new matrices of the isometries which send horizontal planes to horizontal planes for $C_{0}$}
 \label{table:motionsC}
 \end{table}


\begin{table}
\small
\[
\begin{array}{|c|c|c|c|c|}
\hline
& \multicolumn{4}{c|}{\text{\normalsize Type {\em D}}} \\  
\hline
\text{\normalsize $\mathcal{C}_{1}$} &
 
\begin{matrix}
1 & 0 & 0 & 0 \\
1 & 1 & 0 & 0 \\
1 & 0 & 1 & 0 \\
1 & 0 & 0 & 1 
\end{matrix}
 
 & 
  
\begin{matrix}
1 & 0 & 0 & 0 \\
0 & -1 & 0 & 0 \\
1/2 & 0 & -1 & 0 \\
1/2 & 0 & 0 & 1 
\end{matrix}
 
 & 
  
\begin{matrix}
1 & 0 & 0 & 0 \\
1/2 & -1 & 0 & 0 \\
1/2 & 0 & 1 & 0 \\
0 & 0 & 0 & -1 
\end{matrix}
 
&
  
\begin{matrix}
1 & 0 & 0 & 0 \\
1/2 & 1 & 0 & 0 \\
0 & 0 & -1 & 0 \\
1/2 & 0 & 0 & -1 
\end{matrix}
 
 \\
 \hline
 %
 %
%
%
%
%
%
%
%
%
\text{\normalsize $\mathcal{C}_{3}$}  &
 
\begin{matrix}
1 & 0 & 0 & 0 \\
0 & -1 & 0 & 0 \\
1 & 0 & 1 & 0 \\
1/2 & 0 & 0 & -1 
\end{matrix}
 
&
  
\begin{matrix}
1 & 0 & 0 & 0 \\
1 & 1 & 0 & 0 \\
1/2 & 0 & -1 & 0 \\
0 & 0 & 0 & -1 
\end{matrix}
 
 & 
%
%
&
 
\begin{matrix}
1 & 0 & 0 & 0 \\
1/2 & -1 & 0 & 0 \\
0 & 0 & -1 & 0 \\
1 & 0 & 0 & 1 
\end{matrix}
 
\\
 \hline
 %
 %
%
%
%
%
%
%
%
%
\text{\normalsize $\mathcal{C}_{5}$}  &
 
\begin{matrix}
1 & 0 & 0 & 0 \\
1/2 & 1 & 0 & 0 \\
1 & 0 & 1 & 0 \\
0 & 0 & 0 & -1 
\end{matrix}
 
& 
%
%
 &
 
\begin{matrix}
1 & 0 & 0 & 0 \\
0 & -1 & 0 & 0 \\
1/2 & 0 & 1 & 0 \\
1 & 0 & 0 & 1 
\end{matrix}
 
 &
   
\begin{matrix}
1 & 0 & 0 & 0 \\
1 & 1 & 0 & 0 \\
0 & 0 & -1 & 0 \\
1/2 & 0 & 0 & 1 
\end{matrix}
 
\\
 \hline
 %
 %
%
%
%
%
%
%
%
%
 \text{\normalsize $\mathcal{C}_{7}$}  &
 
\begin{matrix}
1 & 0 & 0 & 0 \\
1/2 & -1 & 0 & 0 \\
1 & 0 & 1 & 0 \\
1/2 & 0 & 0 & 1 
\end{matrix}
 
 &
  
\begin{matrix}
1 & 0 & 0 & 0 \\
1/2 & 1 & 0 & 0 \\
1/2 & 0 & -1 & 0 \\
1 & 0 & 0 & 1 
\end{matrix}
 
& 
  
\begin{matrix}
1 & 0 & 0 & 0 \\
1 & 1 & 0 & 0 \\
1/2 & 0 & 1 & 0 \\
1/2 & 0 & 0 & -1 
\end{matrix}
 
 & 
 
\begin{matrix}
1 & 0 & 0 & 0 \\
0 & -1 & 0 & 0 \\
0 & 0 & -1 & 0 \\
0 & 0 & 0 & -1 
\end{matrix}
 
 \\
 \hline
 %
 %
%
%
%
%
%
%
%
%
\end{array}
 \]
 \caption{The new matrices of the isometries which send horizontal planes to horizontal planes for $D_{0}$}
 \label{table:motionsD}
 \end{table}

\subsection{Cutting with translations and rotations}

At this point we can say that we have already computed a single extended Voronoi region for the eight groups and the four prototiles: the cube $[-1,1]^3$ decomposed into $768\times 4= 3072$ tiles, 768 of each type $A$, $B$, $C$ or $D$. 
But this region is certainly not good enough for our purposes.
%
%
Our next step is to apply Lemmas~\ref{lemma:VorExt-rotation} and~\ref{lemma:VorExt-translation} to
exclude from this initial population as many tiles as we can.

That is to say: Let $G$ be one of the  eight quarter groups  and let $T_0$ be one of the four prototiles $A_0$, $B_0$, $C_0$ and $D_0$. For each choice of $G$ and $T_0$ we do the following:

\begin{itemize}
\item Let $S$ be a set of rotations and translations present in $G$. Below we specify our choice of $S$ for each $G$, but essentially what we do is we take the translations of smaller length and the rotations of axis closer to the origin. 

\item For each of the $768\times 4$ tiles in the initial population and for each element $\rho\in S$, apply 
Lemma~\ref{lemma:VorExt-rotation} (if $\rho$ is a rotation) or~\ref{lemma:VorExt-translation} (if $\rho$ is a translation) and, if the tile is fully contained in the region of space forbidden by the Lemma, then remove this tile from the extended Voronoi region.

There is a subtle point in this latter test. Let $T$ be a tile to which we can apply the test of Lemma, say,~\ref{lemma:VorExt-rotation}. The ``forbidden region'' is the unbounded dihedron defined by two half-spaces $H'_1$ and $H'_2$ (see Figure~\ref{fig:VorExtRotTrans}). What we test is whether each of the vertices of $T$ belongs to each of the two half-spaces. But the decision of whether the tile is discarded or not depends on whether the dihedron is convex or not. That is, whether the angle $\alpha$ of the rotation $\rho$ plus the angle with which $T_0$ is seen from $l$ is bigger or smaller than $180$ degrees:

\begin{itemize}
\item If this angle is bigger than 180 degrees (that is, if the forbidden region is convex) we discard the tile $T$ whenever ``both half-spaces $H'_1$ and $H'_2$ contain all the vertices of $T$".

\item If this angle is smaller than 180 degrees (that is, if the forbidden region is not convex) we discard the tile $T$ whenever ``one of the half-spaces $H'_1$ and $H'_2$ contains all the vertices of $T$".

\end{itemize}

If $\rho$ is a translation, the second rule is always applied.
\end{itemize}

The particular translations and rotations that we use are the following ones. Tables~\ref{table:transformations}, \ref{table:axes-triad-rotations}, \ref{table:coordinate-diad-rotations}, \ref{table:diagonal-diad-rotations} and~\ref{table:diagonal-P4132} list them more explicitly.


\begin{table}
\[
\begin{array}{|c|c|}
 \hline
\text{Sets} & \text{Isommetries}\\
\hline
S_{1} & \begin{array}{c}
             \text{Coordinate translations: } \\ (\pm1,0,0),  (0,\pm1,0),  (0,0,\pm1) \end{array}  \\
\hline
S_{2} &\begin{array}{c}
\text{Diagonal translations: }\\ (\pm1/2,\pm1/2,\pm1/2)  \end{array}
\\
\hline
S_{3} & \text{Triad  rotations  (Table~\ref{table:axes-triad-rotations})} \\
\hline
S_{4} & \text{Diad   rotations  (Table~\ref{table:coordinate-diad-rotations}) }\\
\hline
S_{5} & \text{Diad rotations (Table~\ref{table:diagonal-diad-rotations})} \\
\hline
S_{6} & \text{Diad   rotations  (Table~\ref{table:diagonal-P4132})} \\
\hline
\end{array}
\qquad
\begin{array}{|c|l|l|}
 \hline
\text{Group} & \text{Transl.}& \text{Rotations}\\
\hline
  \begin{array}{c}
  P2_{1}3\\ P\frac{2_{1}}{a}\overline{3} 
  \end{array}
  & S_{1} &S_{3} \\
\hline
P4_{1}32 & S_{1} &S_{3}, S_{6}\\
\hline
\begin{array}{c}
I\overline{4}3d\\ I\frac{2}{g}\overline{3}\\ I2'3 
\end{array}
& S_{1}, S_{2} &S_{3}, S_{4} \\
\hline
\begin{array}{c}
I\frac{4_{1}}{g}\overline{3}\frac{2}{d}\\ I4_{1}32 
\end{array}
& S_{1}, S_{2} &S_{3}, S_{4}, S_{5} \\
\hline
\end{array}
\]
\caption{Transformations that we use in each quarter group}
\label{table:transformations}
\end{table}   


\begin{itemize}

\item The integer translations of length one appear in all the quarter groups. In the five groups with a body centered translational subgroup we also use the half-integer diagonal translations $(\pm1/2,\pm1/2,\pm1/2)$. 

\item The triad rotations also appear in all the quarter groups. The set $S_3$ of  triad rotations that we use are listed in Table~\ref{table:axes-triad-rotations}. We have included the ones with axis closer to the prototile, since the ones that are not close
will cut out only portions of the extended Voronoi region that were already cut by other, closer rotations.
The decision on where to put the threshold is somehow heuristic but should not affect the result much. Actually, many of the rotations listed in our table turned out to be already superfluous.


\begin{table}
\[
\begin{array}{cc|cc|cc}
\multicolumn{6}{c}{\text{\normalsize Triad rotations}} \\
\hline
\text{Point} & \text{Vector} & \text{Point} & \text{Vector} & \text{Point} & \text{Vector}\\
\hline
(0,0,0) & (1,1,1) & (-1/2,0,0) & (-1,-1,1) & (0,1/2,0) & (-1,1,1) \\
(-1,0,0) & (1,1,1) & (1/2,0,0) & (-1,-1,1) & (0,-1/2,0) & (-1,1,1) \\
(1,0,0) & (1,1,1)  & (1/2,1,0) & (-1,-1,1) & (-1,1/2,0) & (-1,1,1) \\
(0,-1,0) & (1,1,1)  & (-1/2,1,0) & (-1,-1,1) & (-1,-1/2,0) & (-1,1,1) \\
(0,1,0) & (1,1,1) & (1/2,-1,0) & (-1,-1,1) & (1,1/2,0) & (-1,1,1) \\
(0,0,-1) & (1,1,1)  & (-1/2,-1,0) & (-1,-1,1) & (1,-1/2,0) & (-1,1,1) \\
(0,0,1) & (1,1,1)  & (-3/2,-1,0) & (-1,-1,1) & (-1,3/2,0) & (-1,1,1) \\
(0,0,1/2) & (1,-1,1)  & (1/2,0,1) & (-1,-1,1) & (0,1/2,-1) & (-1,1,1) \\
(1,0,1/2) & (1,-1,1) & (-1/2,0,1) & (-1,-1,1) & (0,-1/2,-1) & (-1,1,1) \\
(0,0,-1/2) & (1,-1,1) & (1/2,0,-1) & (-1,-1,1) & (0,1/2,1) & (-1,1,1) \\
(1,0,-1/2) & (1,-1,1) & (-1/2,0,-1) & (-1,-1,1) & (0,-1/2,1) & (-1,1,1) \\
(0,1,-1/2) & (1,-1,1)  & (0,0,3/2) & (1,-1,1) &  (1,-1/2,0) & (-1,1,1)\\
(0,1,1/2) & (1,-1,1) & (0,0,-3/2) & (1,-1,1) & (1,1/2,0) & (-1,1,1) \\
 &  & (1/2,3/2,0) & (1,-1,1) & (-1,1/2,0) & (-1,1,1) \\
%
\end{array}
\]
\caption{A set $S_3$ of triad rotations in $Q$ and, therefore, in all the quarter groups}
\label{table:axes-triad-rotations}
\end{table} 

\begin{table}
\[
\begin{array}{cc|cc|cc}
\multicolumn{6}{c}{\text{\normalsize Diad rotations parallel to the coordinate axes}}\\
\hline
\text{Point} & \text{Vector} & \text{Point} & \text{Vector} & \text{Point} & \text{Vector}\\
\hline
(0,-3/4,-1/2) & (1,0,0) & (-1/2,0,-3/4) & (0,1,0) & (-3/4,-1/2,0) & (0,0,1)\\
(0,-3/4,0) & (1,0,0) & (-1/2,0,-1/4) & (0,1,0) & (-3/4,0,0) & (0,0,1) \\
(0,-3/4,1/2) & (1,0,0) & (-1/2,0,1/4) & (0,1,0) & (-3/4,1/2,0) & (0,0,1) \\
(0,-1/4,-1/2) & (1,0,0) & (-1/2,0,3/4) & (0,1,0) & (-1/4,-1/2,0) & (0,0,1)\\
(0,-1/4,0) & (1,0,0) &  (0,0,-3/4) & (0,1,0) & (-1/4,0,0) & (0,0,1) \\
 (0,-1/4,1/2) & (1,0,0) & (0,0,-1/4) & (0,1,0) & (-1/4,1/2,0) & (0,0,1)\\
(0,1/4,-1/2) & (1,0,0) & (0,0,1/4) & (0,1,0) & (1/4,-1/2,0) & (0,0,1)\\
 (0,1/4,0) & (1,0,0) & (0,0,3/4) & (0,1,0) & (1/4,0,0) & (0,0,1)\\
 (0,1/4,1/2) & (1,0,0) & (1/2,0,-3/4) & (0,1,0) & (1/4,1/2,0) & (0,0,1)\\
(0,3/4,-1/2) & (1,0,0) &  (1/2,0,-1/4) & (0,1,0) & (3/4,-1/2,0) & (0,0,1) \\
 (0,3/4,0) & (1,0,0) & (1/2,0,1/4) & (0,1,0) & (3/4,0,0) & (0,0,1) \\
(0,3/4,1/2) & (1,0,0) & (1/2,0,3/4) & (0,1,0) & (3/4,1/2,0) & (0,0,1) \\
\end{array}
\]
\caption{A set $S_4$ of diad rotations parallel to the coordinate axes that appear in the groups $I\frac{4_{1}}{g}\overline{3}\frac{2}{d}$, $I4_{1}32$, $I43d$, $I\frac{2}{g}\overline{3}$ and $I2'3$}
\label{table:coordinate-diad-rotations}
\end{table}  

\begin{table}
\footnotesize
\[
\begin{array}{cc|cc|cc}
\hline
\text{Point} & \text{Vector} & \text{Point} & \text{Vector} & \text{Point} & \text{Vector}\\
\hline
(-7/4,0,1/8) & (-1,1,0) & (-5/4,0,3/8) & (1,1,0) & (1/8,-7/4,0) & (0,-1,1)\\
(-7/4,0,5/8) & (-1,1,0) & (-5/4,0,7/8) & (1,1,0) & (5/8,-7/4,0) & (0,-1,1) \\
(-5/4,0,3/8) & (-1,1,0) & (-3/4,0,1/8) & (1,1,0) & (3/8,-5/4,0) & (0,-1,1) \\
(-5/4,0,7/8) & (-1,1,0) & (-3/4,0,5/8) & (1,1,0) & (7/8,-5/4,0) & (0,-1,1)\\
(-3/4,0,1/8) & (-1,1,0) &  (-1/4,0,3/8) & (1,1,0) & (1/8,-3/4,0) & (0,-1,1) \\
 (-3/4,0,5/8) & (-1,1,0) & (-1/4,0,7/8) & (1,1,0) & (5/8,-3/4,0) & (0,-1,1)\\
(-1/4,0,3/8) & (-1,1,0) & (1/4,0,1/8) & (1,1,0) & (3/8,-1/4,0) & (0,-1,1)\\
 (-1/4,0,7/8) & (-1,1,0) & (1/4,0,5/8) & (1,1,0) & (7/8,-1/4,0) & (0,-1,1)\\
 (1/4,0,1/8) & (-1,1,0) & (3/4,0,3/8) & (1,1,0) & (1/8,1/4,0) & (0,-1,1)\\
(1/4,0,5/8) & (-1,1,0) &  (3/4,0,7/8) & (1,1,0) & (5/8,1/4,0) & (0,-1,1) \\
 (3/4,0,3/8) & (-1,1,0) & (5/4,0,1/8) & (1,1,0) & (3/8,3/4,0) & (0,-1,1) \\
(3/4,0,7/8) & (-1,1,0) & (5/4,0,5/8) & (1,1,0) & (7/8,3/4,0) & (0,-1,1) \\
(5/4,0,1/8) & (-1,1,0) & (7/4,0,3/8) & (1,1,0) & (1/8,5/4,0) & (0,-1,1) \\
(5/4,0,5/8) & (-1,1,0) & (7/4,0,7/8) & (1,1,0) & (5/8,5/4,0) & (0,-1,1) \\
(7/4,0,3/8) & (-1,1,0) & (-5/4,0,-1/8) & (1,1,0) & (3/8,7/4,0) & (0,-1,1) \\
(7/4,0,7/8) & (-1,1,0) &(-5/4,0,-5/8) & (1,1,0)& (7/8,7/4,0) & (0,-1,1) \\
(-7/4,0,-3/8) & (-1,1,0) & (-3/4,0,-3/8) & (1,1,0)  & (-3/8,-7/4,0) & (0,-1,1)\\
(-7/4,0,-7/8) & (-1,1,0) &(-3/4,0,-7/8) & (1,1,0) & (-7/8,-7/4,0) & (0,-1,1) \\
(-5/4,0,-1/8) & (-1,1,0) & (-1/4,0,-1/8) & (1,1,0) & (-1/8,-5/4,0) & (0,-1,1) \\
(-5/4,0,-5/8) & (-1,1,0) & (-1/4,0,-5/8) & (1,1,0) & (-5/8,-5/4,0) & (0,-1,1)\\
(-3/4,0,-3/8) & (-1,1,0) &  (1/4,0,-3/8) & (1,1,0) & (-3/8,-3/4,0) & (0,-1,1) \\
 (-3/4,0,-7/8) & (-1,1,0) &  (1/4,0,-7/8) & (1,1,0) & (-7/8,-3/4,0) & (0,-1,1)\\
(-1/4,0,-1/8) & (-1,1,0) & (3/4,0,-1/8) & (1,1,0) & (-1/8,-1/4,0) & (0,-1,1)\\
 (-1/4,0,-5/8) & (-1,1,0) & (3/4,0,-5/8) & (1,1,0) & (-5/8,-1/4,0) & (0,-1,1)\\
 (1/4,0,-3/8) & (-1,1,0) & (5/4,0,-3/8) & (1,1,0) & (-3/8,1/4,0) & (0,-1,1)\\
(1/4,0,-7/8) & (-1,1,0) & (5/4,0,-7/8) & (1,1,0)  & (-7/8,1/4,0) & (0,-1,1) \\
 (3/4,0,-1/8) & (-1,1,0) & (7/4,0,-1/8) & (1,1,0)  & (-1/8,3/4,0) & (0,-1,1) \\
(3/4,0,-5/8) & (-1,1,0) & (7/4,0,-5/8) & (1,1,0) & (-5/8,3/4,0) & (0,-1,1) \\
(5/4,0,-3/8) & (-1,1,0) & (0,1/8,-7/4) & (1,0,-1)  & (-3/8,5/4,0) & (0,-1,1) \\
(5/4,0,-7/8) & (-1,1,0) & (0,1/8,-7/4) & (1,0,-1)  & (-7/8,5/4,0) & (0,-1,1) \\
(7/4,0,-1/8) & (-1,1,0) &   (0,3/8,-5/4) & (1,0,-1)  & (-1/8,7/4,0) & (0,-1,1) \\
(7/4,0,-5/8) & (-1,1,0) &   (0,7/8,-5/4) & (1,0,-1)  & (-5/8,7/4,0) & (0,-1,1) \\
(3/8,-5/4,0) & (0,1,1) &  (0,1/8,-3/4) & (1,0,-1) & (0,3/8,-5/4) & (1,0,1) \\
(7/8,-5/4,0) & (0,1,1) &  (0,5/8,-3/4) & (1,0,-1) & (0,7/8,-5/4) & (1,0,1) \\
(1/8,-3/4,0) & (0,1,1) &  (0,3/8,-1/4) & (1,0,-1)& (0,1/8,-3/4) & (1,0,1) \\
(5/8,-3/4,0) & (0,1,1) &  (0,7/8,-1/4) & (1,0,-1) & (0,5/8,-3/4) & (1,0,1) \\
(3/8,-1/4,0) & (0,1,1) &  (0,1/8,1/4) & (1,0,-1)& (0,3/8,-1/4) & (1,0,1) \\
(7/8,-1/4,0) & (0,1,1) & (0,5/8,1/4) & (1,0,-1) & (0,7/8,-1/4) & (1,0,1) \\
(1/8,1/4,0) & (0,1,1) & (0,3/8,3/4) & (1,0,-1)  & (0,1/8,1/4) & (1,0,1) \\
(5/8,1/4,0) & (0,1,1) & (0,7/8,3/4) & (1,0,-1)  & (0,5/8,1/4) & (1,0,1) \\
(3/8,3/4,0) & (0,1,1) & (0,1/8,5/4) & (1,0,-1) & (0,3/8,3/4) & (1,0,1) \\
(7/8,3/4,0) & (0,1,1) & (0,5/8,5/4) & (1,0,-1)& (0,7/8,3/4) & (1,0,1) \\
(1/8,5/4,0) & (0,1,1) & (0,3/8,7/4) & (1,0,-1) & (0,1/8,5/4) & (1,0,1) \\
(5/8,5/4,0) & (0,1,1) & (0,7/8,7/4) & (1,0,-1) & (0,5/8,5/4) & (1,0,1) \\
(3/8,7/4,0) & (0,1,1) & (0,-3/8,-7/4) & (1,0,-1) & (0,3/8,7/4) & (1,0,1) \\
(7/8,7/4,0) & (0,1,1) & (0,-7/8,-7/4) & (1,0,-1)  & (0,7/8,7/4) & (1,0,1) \\
(-1/8,-5/4,0) & (0,1,1) &  (0,-1/8,-5/4,0) & (1,0,-1) & (0,-1/8,-5/4) & (1,0,1) \\
(-5/8,-5/4,0) & (0,1,1) & (0,-5/8,-5/4) & (1,0,-1) & (0,-5/8,-5/4) & (1,0,1) \\
(-3/8,-3/4,0) & (0,1,1) & (0,-3/8,-3/4) & (1,0,-1) & (0,-3/8,-3/4) & (1,0,1) \\
(-7/8,-3/4,0) & (0,1,1) & (0,-7/8,-3/4) & (1,0,-1)  & (0,-7/8,-3/4) & (1,0,1) \\
(-1/8,-1/4,0) & (0,1,1) & (0,-1/8,-1/4) & (1,0,-1)  & (0,-1/8,-1/4) & (1,0,1) \\
(-5/8,-1/4,0) & (0,1,1) &  (0,-5/8,-1/4) & (1,0,-1) & (0,-5/8,-1/4) & (1,0,1) \\
(-3/8,1/4,0) & (0,1,1) & (0,-3/8,1/4) & (1,0,-1) & (0,-3/8,1/4) & (1,0,1) \\
(-7/8,1/4,0) & (0,1,1) & (0,-7/8,1/4) & (1,0,-1)  & (0,-7/8,1/4) & (1,0,1) \\
(-1/8,3/4,0) & (0,1,1) &  (0,-1/8,3/4) & (1,0,-1) & (0,-1/8,3/4) & (1,0,1) \\
(-5/8,3/4,0) & (0,1,1) & (0,-5/8,3/4) & (1,0,-1) & (0,-5/8,3/4) & (1,0,1) \\
(-3/8,5/4,0) & (0,1,1) & (0,-3/8,5/4) & (1,0,-1) & (0,-3/8,5/4) & (1,0,1) \\
(-7/8,5/4,0) & (0,1,1) & (0,-7/8,5/4) & (1,0,-1) & (0,-7/8,5/4) & (1,0,1) \\
(-1/8,7/4,0) & (0,1,1) & (0,-1/8,7/4) & (1,0,-1) & (0,-1/8,7/4) & (1,0,1) \\
(-5/8,7/4,0) & (0,1,1) & (0,-5/8,7/4) & (1,0,-1)  & (0,-5/8,7/4) & (1,0,1) \\
\end{array}
\]
\caption{A set $S_5$ of diad rotations parallel to the diagonal of the faces of the cube, in the groups $I\frac{4_{1}}{g}\overline{3}\frac{2}{d}$ and $I4_{1}32$}
\label{table:diagonal-diad-rotations}
\end{table}


\begin{table}
\footnotesize
\[
\begin{array}{cc|cc|cc}
\multicolumn{6}{c}{\text{\normalsize Diad rotations parallel to the diagonals of the faces of the unit cube}}\\
\hline
\text{Point} & \text{Vector} & \text{Point} & \text{Vector} & \text{Point} & \text{Vector}\\
\hline
(-5/4,0,3/8) & (-1,1,0) & (-3/4,0,1/8) & (1,1,0) & (3/8,-5/4,0) & (0,-1,1) \\
(-5/4,0,7/8) & (-1,1,0) & (-3/4,0,5/8) & (1,1,0) & (7/8,-5/4,0) & (0,-1,1)\\
(-1/4,0,3/8) & (-1,1,0) & (1/4,0,1/8) & (1,1,0) & (3/8,-1/4,0) & (0,-1,1)\\
 (-1/4,0,7/8) & (-1,1,0) & (1/4,0,5/8) & (1,1,0) & (7/8,-1/4,0) & (0,-1,1)\\
 (3/4,0,3/8) & (-1,1,0) & (5/4,0,1/8) & (1,1,0) & (3/8,3/4,0) & (0,-1,1) \\
(3/4,0,7/8) & (-1,1,0) & (5/4,0,5/8) & (1,1,0) & (7/8,3/4,0) & (0,-1,1) \\
(7/4,0,3/8) & (-1,1,0) & (-3/4,0,-3/8) & (1,1,0) & (3/8,7/4,0) & (0,-1,1) \\
(7/4,0,7/8) & (-1,1,0) & (-3/4,0,-7/8) & (1,1,0) & (7/8,7/4,0) & (0,-1,1) \\
(-5/4,0,-1/8) & (-1,1,0) & (1/4,0,-3/8) & (1,1,0) & (-1/8,-5/4,0) & (0,-1,1) \\
(-5/4,0,-5/8) & (-1,1,0) & (1/4,0,-7/8) & (1,1,0) & (-5/8,-5/4,0) & (0,-1,1) \\
 (-1/4,0,-1/8) & (-1,1,0) & (5/4,0,-3/8) & (1,1,0) &  (-1/8,-1/4,0) & (0,-1,1)\\
 (-1/4,0,-5/8) & (-1,1,0) & (5/4,0,-7/8) & (1,1,0) & (-5/8,-1/4,0) & (0,-1,1)\\
(3/4,0,-1/8) & (-1,1,0) &   (0,3/8,-5/4) & (1,0,-1) & (-1/8,3/4,0) & (0,-1,1) \\
 (3/4,0,-5/8) & (-1,1,0) &  (0,7/8,-5/4) & (1,0,-1) & (-5/8,3/4,0) & (0,-1,1) \\
 (7/4,0,-1/8) & (-1,1,0) & (0,3/8,-1/4) & (1,0,-1) & (-1/8,7/4,0) & (0,-1,1)\\
 (7/4,0,-5/8) & (-1,1,0)& (0,7/8,-1/4) & (1,0,-1) & (-5/8,7/4,0) & (0,-1,1)\\
(1/8,-3/4,0) & (0,1,1) & (0,3/8,3/4) & (1,0,-1) & (0,1/8,-3/4) & (1,0,1) \\
(5/8,-3/4,0) & (0,1,1) & (0,7/8,3/4) & (1,0,-1)  & (0,5/8,-3/4) & (1,0,1) \\
(1/8,1/4,0) & (0,1,1) & (0,3/8,7/4) & (1,0,-1)  & (0,1/8,1/4) & (1,0,1) \\
(5/8,1/4,0) & (0,1,1) & (0,7/8,7/4) & (1,0,-1)  & (0,5/8,1/4) & (1,0,1) \\
(1/8,5/4,0) & (0,1,1) & (0,-1/8,-5/4,0) & (1,0,-1) & (0,1/8,5/4) & (1,0,1) \\
(5/8,5/4,0) & (0,1,1) & (0,-5/8,-5/4) & (1,0,-1) & (0,5/8,5/4) & (1,0,1) \\
(-3/8,-3/4,0) & (0,1,1) & (0,-1/8,3/4) & (1,0,-1)  & (0,-3/8,-3/4) & (1,0,1) \\
(-7/8,-3/4,0) & (0,1,1) & (0,-5/8,3/4) & (1,0,-1) & (0,-7/8,-3/4) & (1,0,1) \\
(-3/8,1/4,0) & (0,1,1) & (0,-1/8,-1/4) & (1,0,-1) & (0,-3/8,1/4) & (1,0,1) \\
(-7/8,1/4,0) & (0,1,1) &  (0,-5/8,-1/4) & (1,0,-1) & (0,-7/8,1/4) & (1,0,1) \\
(-3/8,5/4,0) & (0,1,1) & (0,-1/8,7/4) & (1,0,-1)  & (0,-3/8,5/4) & (1,0,1) \\
(-7/8,5/4,0) & (0,1,1) & (0,-5/8,7/4) & (1,0,-1) & (0,-7/8,5/4) & (1,0,1) \\
\end{array}
\]
\caption{A set $S_6\subset S_5$ of diad rotations parallel to the diagonal of the faces of the unit cube that appear in the group $P4_{1}32$}
\label{table:diagonal-P4132}
\end{table}

\item Table~\ref{table:coordinate-diad-rotations} lists the set $S_4$ of diad rotations with axes parallel to the coordinate axes that we use in the groups $I\frac{4_{1}}{g}\overline{3}\frac{2}{d}$, $I4_{1}32$, $I\overline{4}3d$, $I\frac{2}{g}\overline{3}$ and $I2'3$.

\item Table~\ref{table:diagonal-diad-rotations} lists the set $S_5$ of diad rotations with axes parallel to the diagonals of the faces of the unit cube that we use in the groups  $I\frac{4_{1}}{g}\overline{3}\frac{2}{d}$ and $I4_{1}32$. In the group $P4_{1}32$ only some of these diad rotations are present. The list $S_6$ of them is in Table~\ref{table:diagonal-P4132}.
\end{itemize}

In these four tables (Tables~\ref{table:axes-triad-rotations}, \ref{table:coordinate-diad-rotations}, \ref{table:diagonal-diad-rotations} and~\ref{table:diagonal-P4132}), the two entries in each column are a point and the direction of the rotation axis.

\subsection{Planar projection}
  \label{sec:final-bounds}
 
 In the extended Voronoi region computed so far only rotations and translations present in our group are taken into account. If we wish to deal with other transformations (glide reflections, screw rotations, etc) we could state the analogues of Lemmas~\ref{lemma:VorExt-rotation} and \ref{lemma:VorExt-translation} for them, but the result would be too complicated, perhaps impractical. 
 
\begin{figure}[htb]
\[
\begin{array}{ccccc}
&&
\includegraphics[width=3cm]{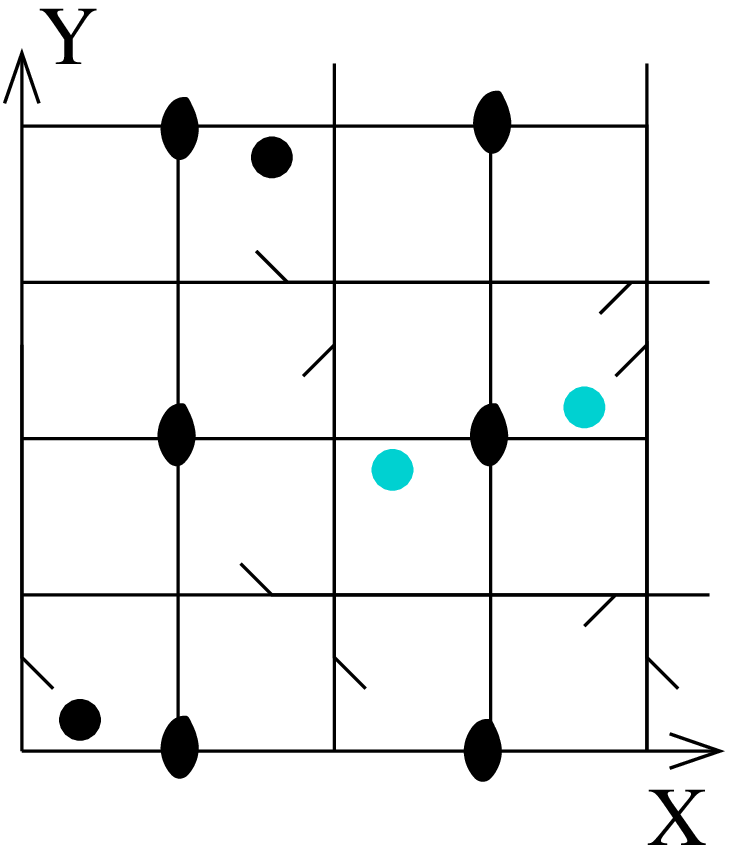}
&\\

&&
pgg
&\\

&\swarrow&\downarrow&\searrow \medskip\\

\includegraphics[width=3cm]{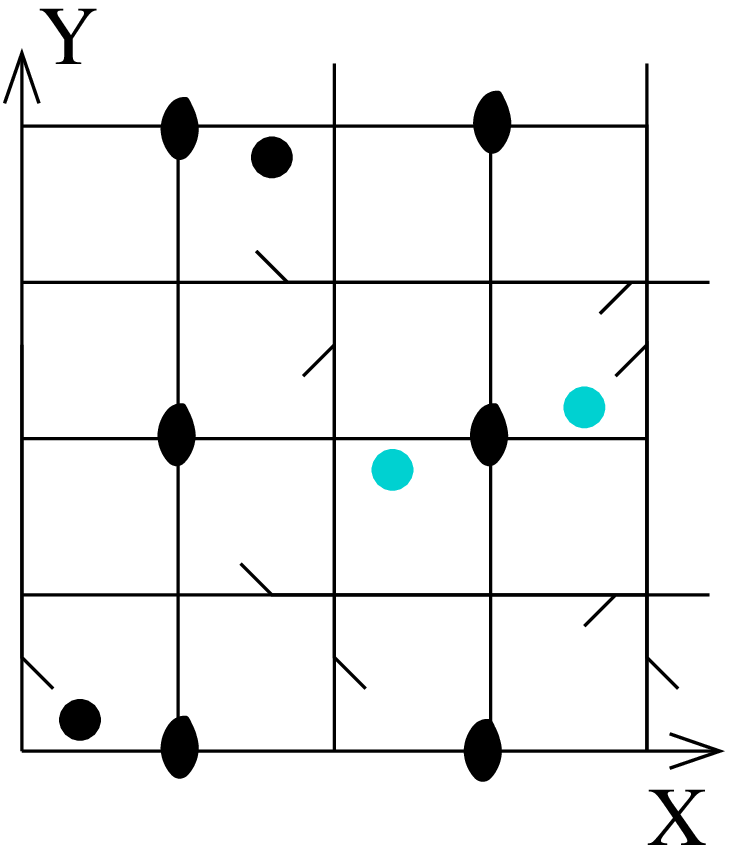}
&&
\includegraphics[width=3cm]{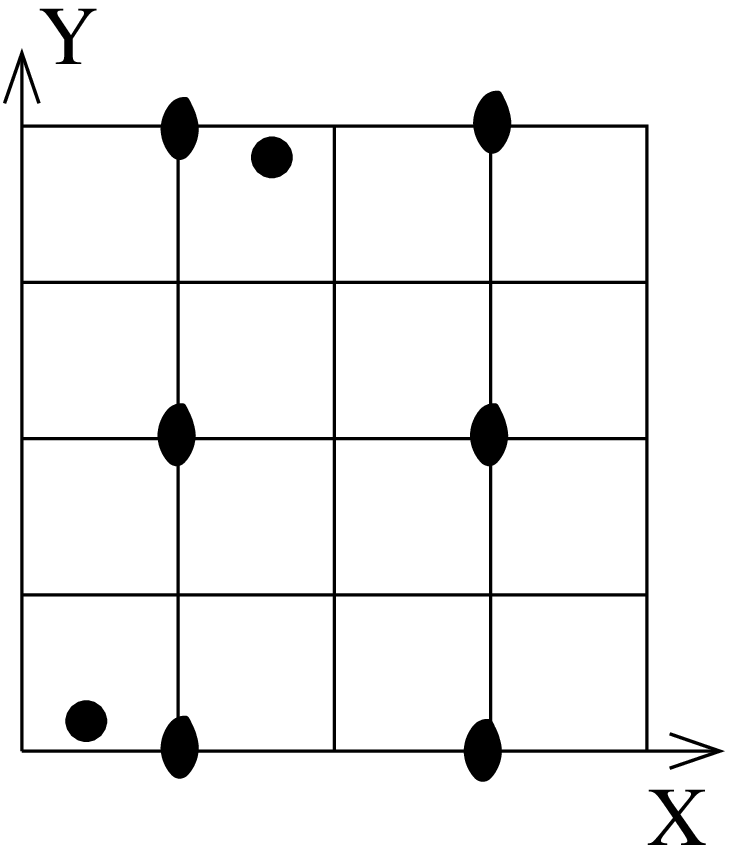}
&&
\includegraphics[width=3cm]{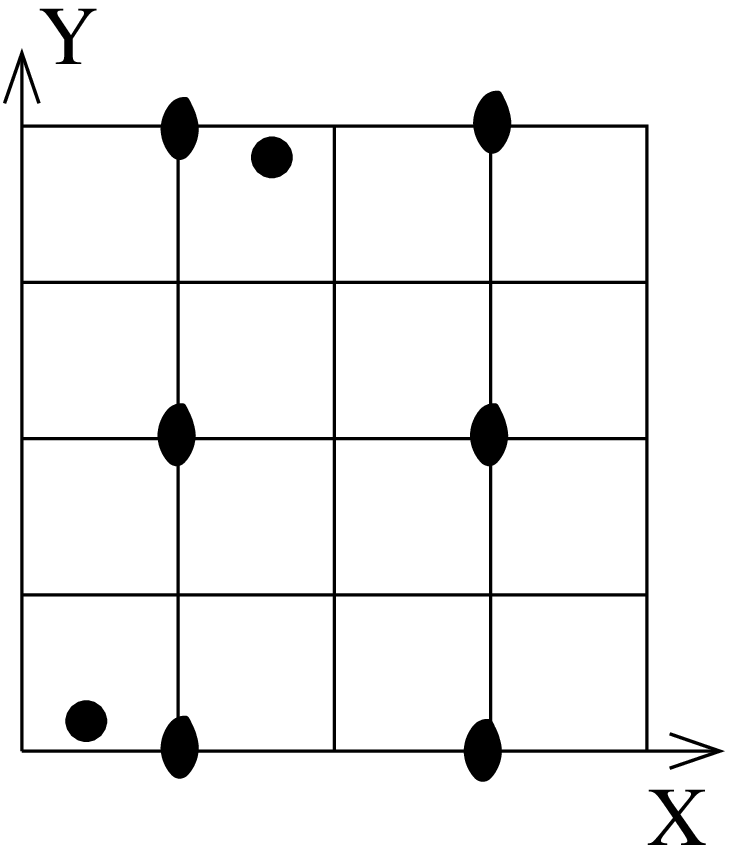}
\\

pgg
&&
p2
&&
p2
\\

\downarrow&\searrow&\downarrow&\swarrow&\downarrow \medskip\\

\includegraphics[width=3cm]{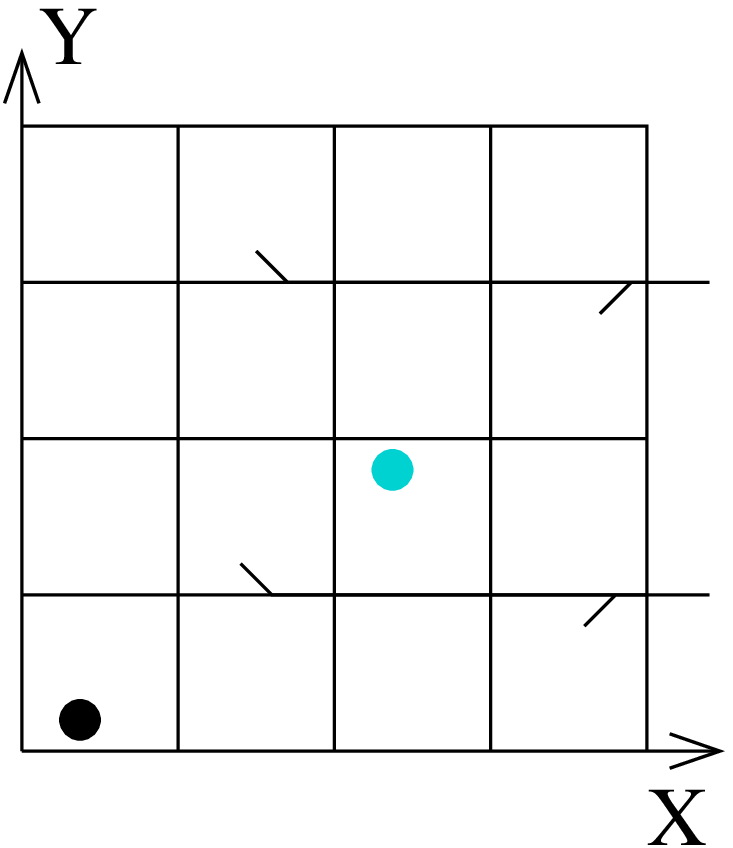}
&&
\includegraphics[width=3cm]{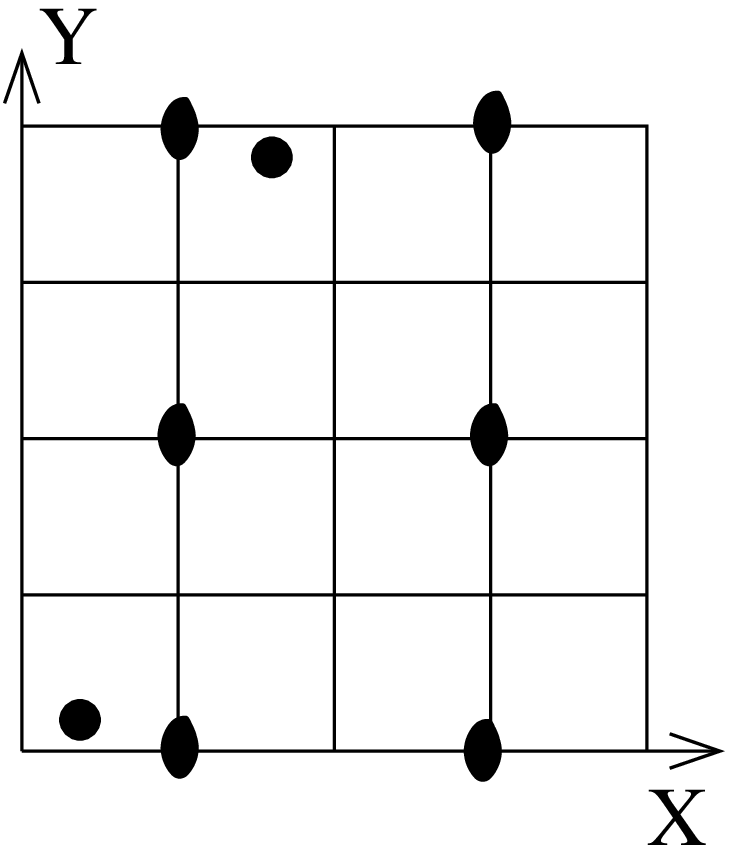}
&&
\includegraphics[width=3cm]{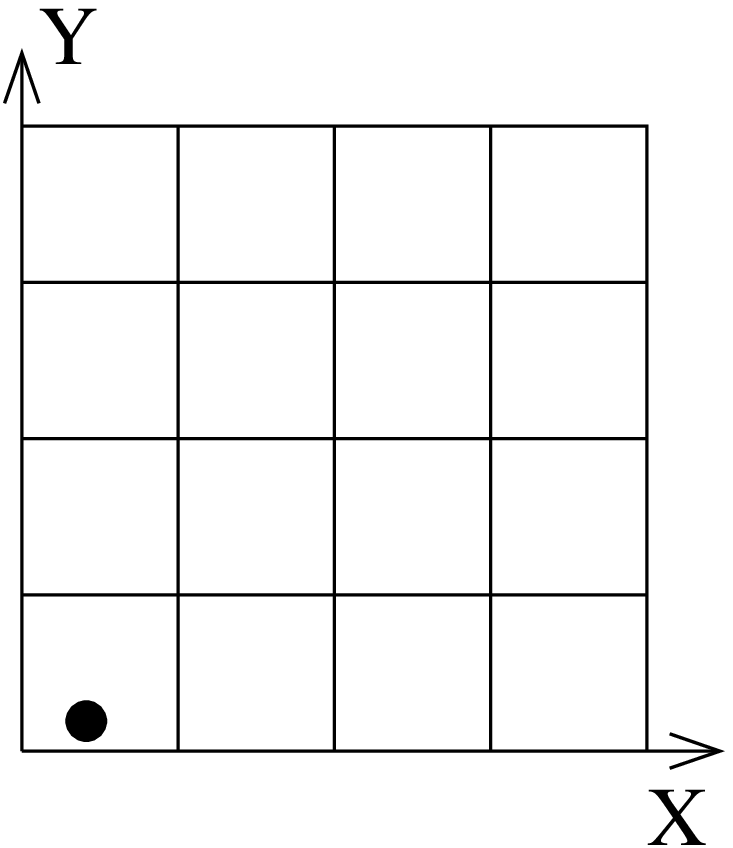}
\\
pg
&&
p2
&&
p1
\\
\medskip

&\searrow&\downarrow&\swarrow&\\

&&
\includegraphics[width=3cm]{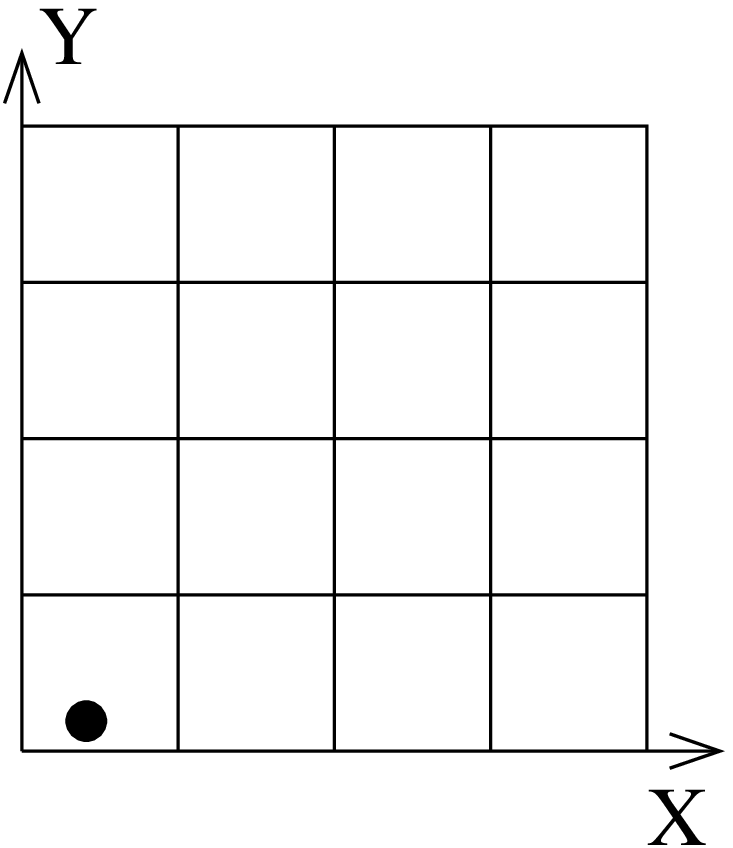}
&&
\\

&&
p1
&&
\\

\end{array}
\]
\caption{Plane groups of the quarter groups.
\label{fig:planegroups}}
\end{figure}

Instead, let us consider the ``horizontal subgroup'' of our group $G$. That is, the subgroup $G_z$ that fixes the third coordinate.  Each of the eight parts of Figure~\ref{fig:planegroups} shows the planar subgroup of the quarter group that occupies the same position in Figure~\ref{fig:symmetries}.
Only four different planar subgroups arise, of types $p1$, $p2$, $pg$ and $pgg$.

As it turns out, extended Voronoi regions (and influence regions) for these four groups were computed in \cite{Bochis-1999} and \cite{Bochis-Santos-2001}. Our idea is to take advantage of the following simple fact:

\begin{lemma}
\label{lemma:project}
Let $\pi:\R^3\to \R^2$ denote the vertical projection, omitting the third coordinate. Let $T_0$ be a prototile and let $G$ be a crystallographic group. Let $G_z$ denote the horizontal subgroup of $G$, that is the set of transformations that send each horizontal plane to itself.
Then:
\[
\pi(\VorExt_{G}(T_0)) \subseteq \VorExt_{G_z}(\pi (T_0)).
\]
\end{lemma}

\begin{proof}
The only thing to prove is the same formula for Voronoi regions. That is, for every base point $p\in T_0$ we need to show that:
\[
\pi(\Vor_{Gp}(p))\subseteq \Vor_{G_z(\pi(p))}(\pi(p)).
\]
The latter follows from the fact that $G_z(\pi(p))$ equals the projection of the part of $Gp$ in the same horizontal plane as $p$ (for a sufficiently generic $p$).
\end{proof}

That is to say: we can intersect our previously computed extended Voronoi regions $\VorExt_{G}(T_0)$
(for $T_0\in \{A_0,\dots,D_0\}$) with $\pi^{-1}\left(\VorExt_{G_z}(\pi (T_0))\right)$ and use the new, smaller region to compute influence regions. In order to do this we check, for each tile in $\VorExt_{G}(T_0)$, whether its projection intersects $\VorExt_{G_z}(\pi (T_0))$. If it does not then we discard that tile. Needless to say that we apply this process three times to each prototile, for the three coordinate projections $XY$, $XZ$ and $YZ$.

For the groups with a planar $p1$ this step does not reduce the regions at all, so we do not show the details. But for the other three types of groups, $p2$, $pg$ and $pgg$, this step is significant.
 Roughly speaking, it allows us to use some glide reflections present in some of the quarter groups to cut the extended Voronoi region further.

\begin{figure}[htb]
\begin{center}
\[\begin{array}{ccc}
\includegraphics[height=4cm]{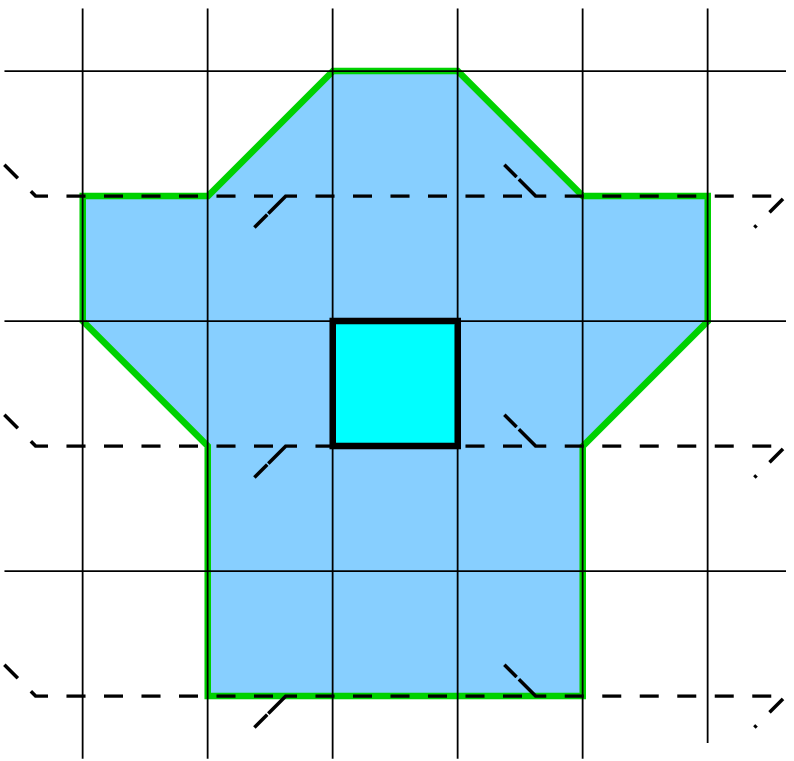}
&
\includegraphics[height =4cm]{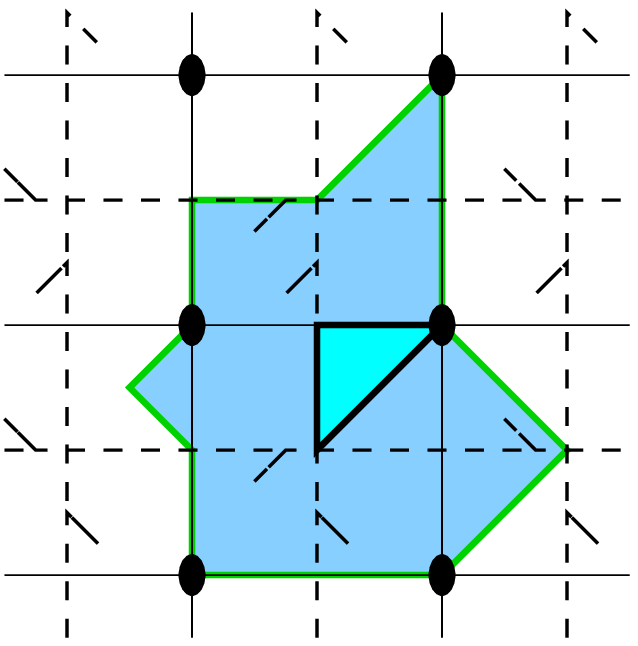} 
&
\includegraphics[height =4cm]{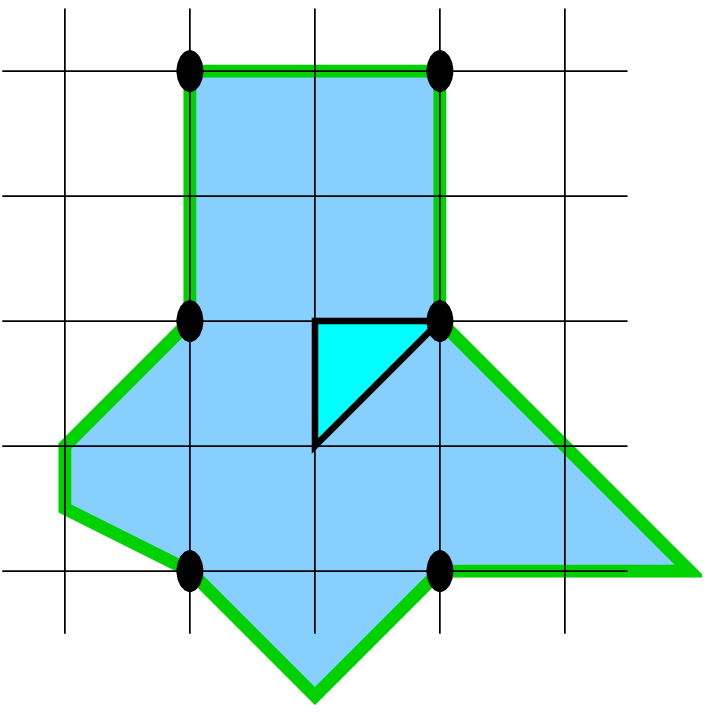}
\\
\text{(a) plane group $pg$.}&
\text{(b) plane group $pgg$.}&
\text{(c) plane group $p2$.}\\
\end{array}
\]
\end{center}
\caption{Extended Voronoi Regions of the plane groups $pg$, $pgg$ and $p2$.
\label{fig:planar-vorext}}
\end{figure}

As a planar prototile we have used a square of side $1/4$ for $pg$ and half of it (a right triangle) for $pgg$ and $p2$. The reason for this difference is that in $pg$ the extended Voronoi region would not decrease significantly if we took a triangle. We omit the calculation of the extended Voronoi regions, shown in Figure~\ref{fig:planar-vorext}. The details can be found in \cite{Bochis-1999} and \cite{Bochis-Santos-2001}. Let us only say that the calculation is based on the planar version of  Lemmas~\ref{lemma:VorExt-rotation} and \ref{lemma:VorExt-translation}, plus the analogue for glide reflections. Roughly speaking, each glide reflection with translation vector, say $v$, ``acts as if'' the group $G_z$ contained a rotation of order two with center on the glide reflection axis and at distance $|v/2|$ from the orthogonal projection of the prototile $T_0$ to the axis. Glide reflections are indicated as dotted lines in the figures.

Of course, we need to check whether each prototile $T_0$ that we used in the 3D step projects to lie in one planar tile. That is, in one tile obtained from the one in Figure~\ref{fig:planar-vorext} by the action of the planar group $G_z$. If this is not the case then the planar extended Voronoi region we consider is going to be a union of more than one of the regions shown in the figure. The result of this check, for each of the four tiles $A_0$, $B_0$, $C_0$ and $D_0$, each of the three coordinate projections, and the three possible plane groups, is illustrated in Tables~\ref{table:intersections} and~\ref{table:intersections2}. As shown in the pictures, for the group $pg$ one planar tile is always enough, but for the groups $pgg$ and $p2$ we usually need two.

 \begin{table}[htb]

\begin{tabular}{|c|c|c|c|}
\hline
 & \multicolumn{3}{c|}{{\text{Type} {\em A}}} \\
\hline
 \footnotesize pgg
&
\includegraphics[width=3cm]{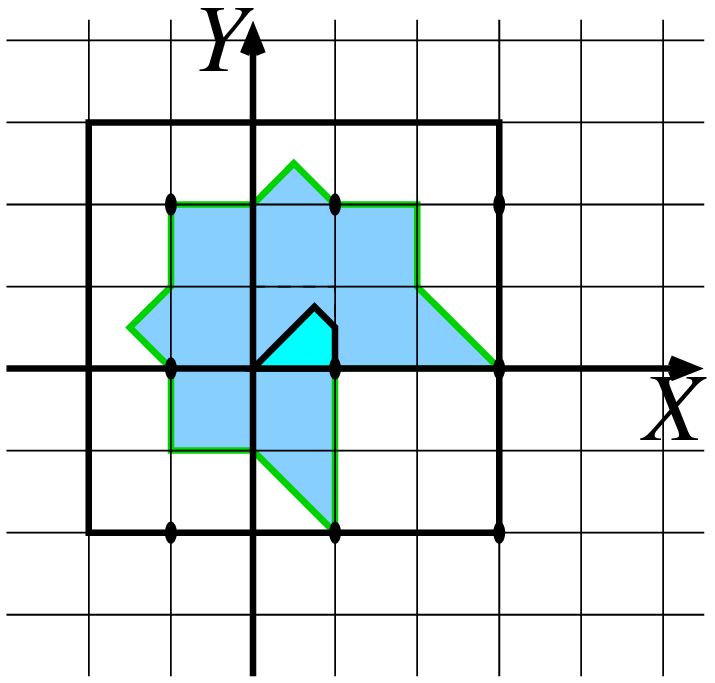}
&
\includegraphics[width=3cm]{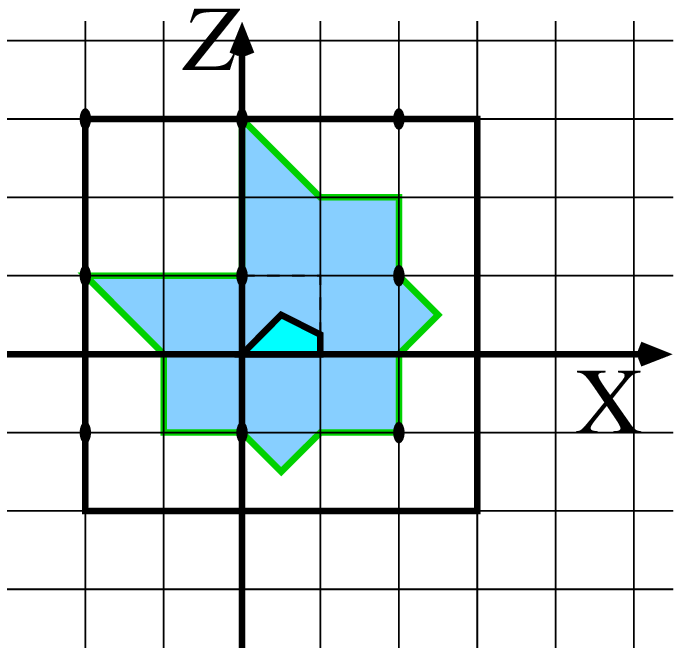}
&
\includegraphics[width=3cm]{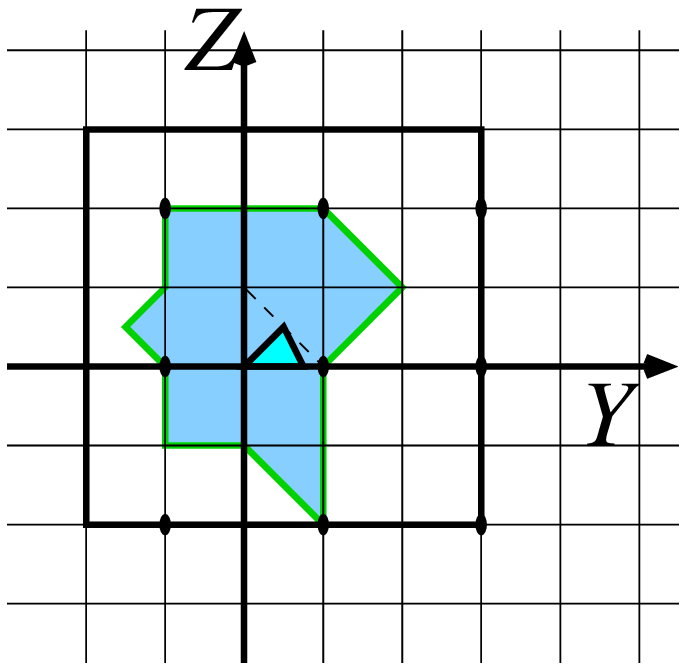}
\\
\hline
\footnotesize p2
&
\includegraphics[width=3cm]{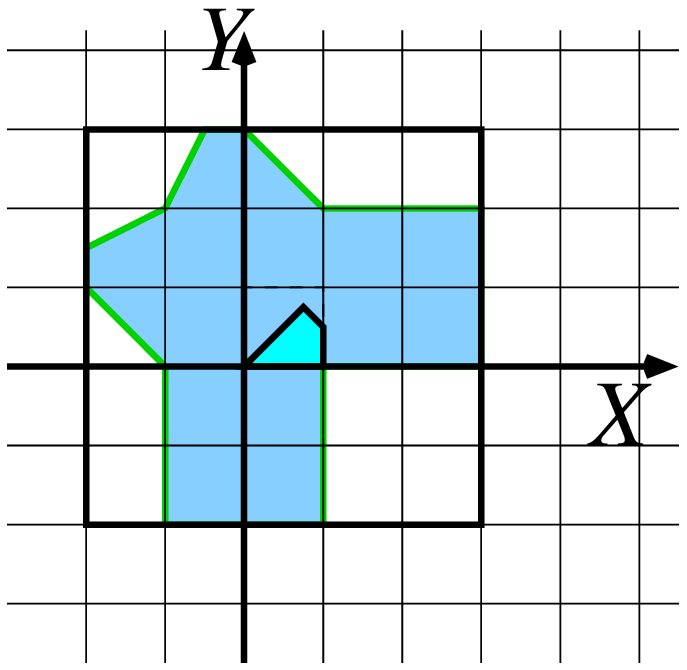}
&
\includegraphics[width=3cm]{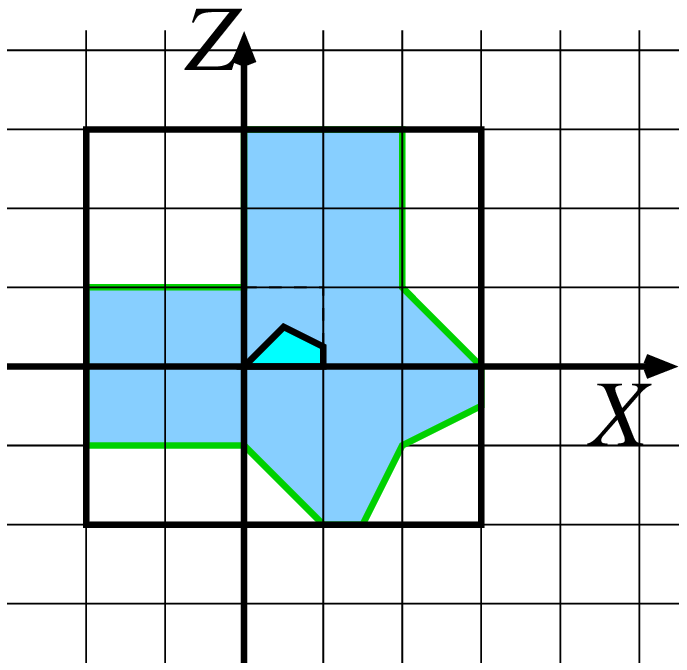} 
&
\includegraphics[width=3cm]{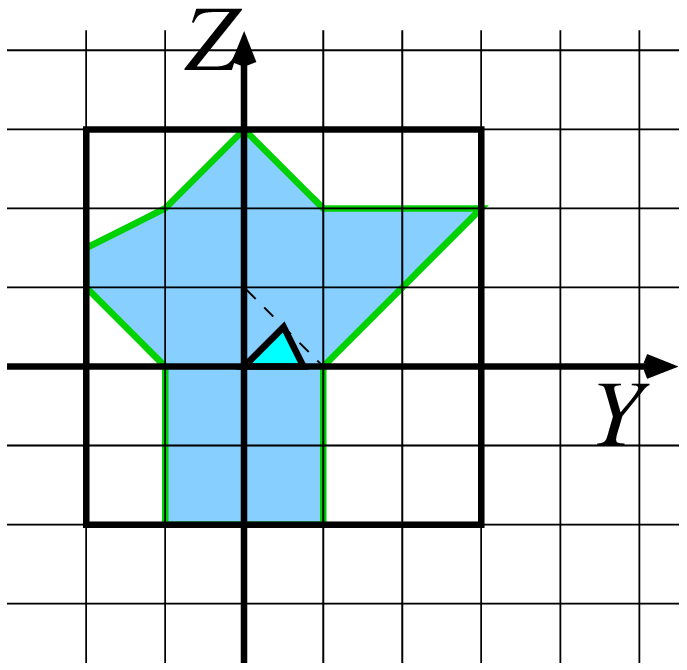} 
\\
\hline
\footnotesize pg
&
\includegraphics[width=3cm]{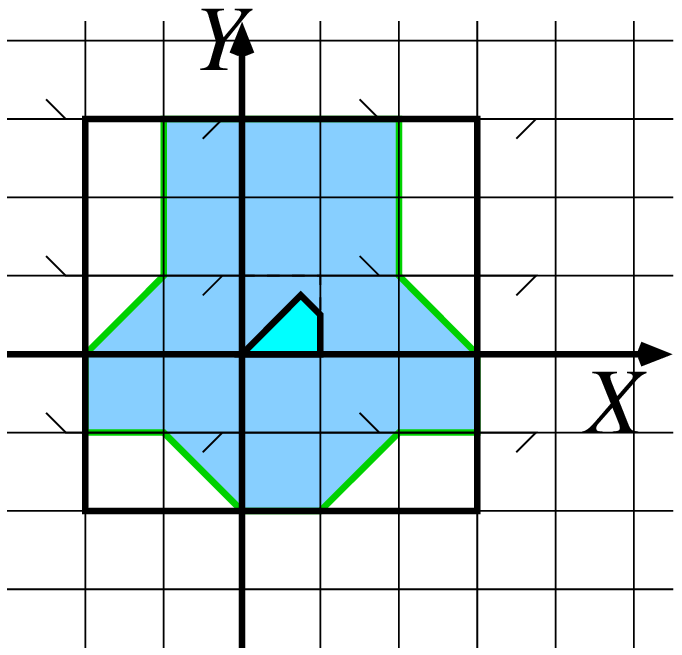}
&
\includegraphics[width=3cm]{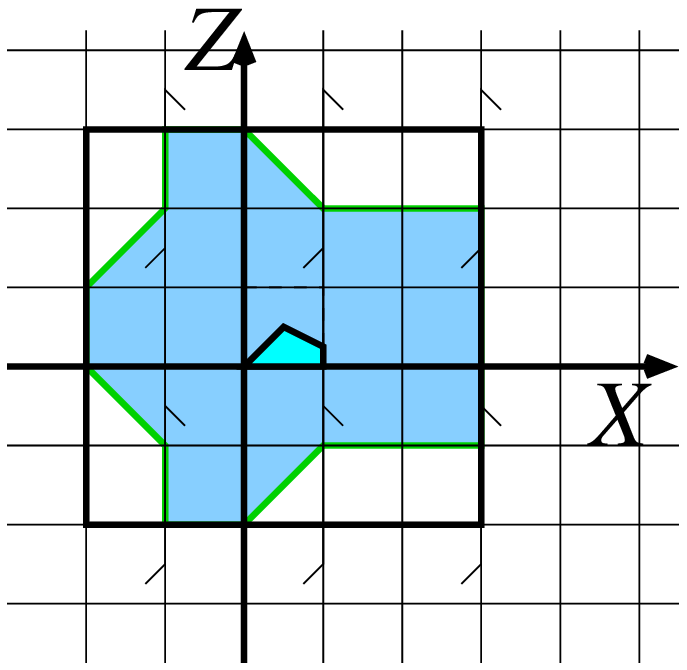}
&
\includegraphics[width=3cm]{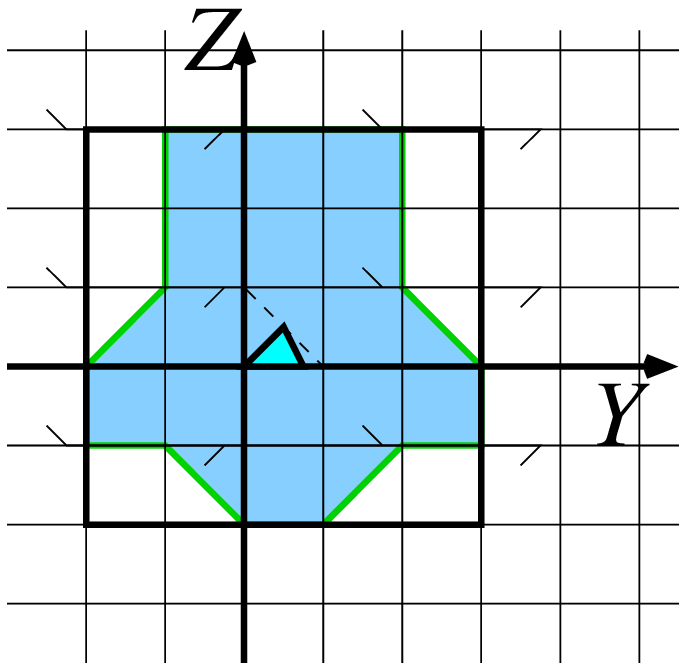}
\\
\hline
\hline
 &  \multicolumn{3}{c|}{{\text{Type} {\em B}}}\\
\hline
 \footnotesize pgg
&
\includegraphics[width=3cm]{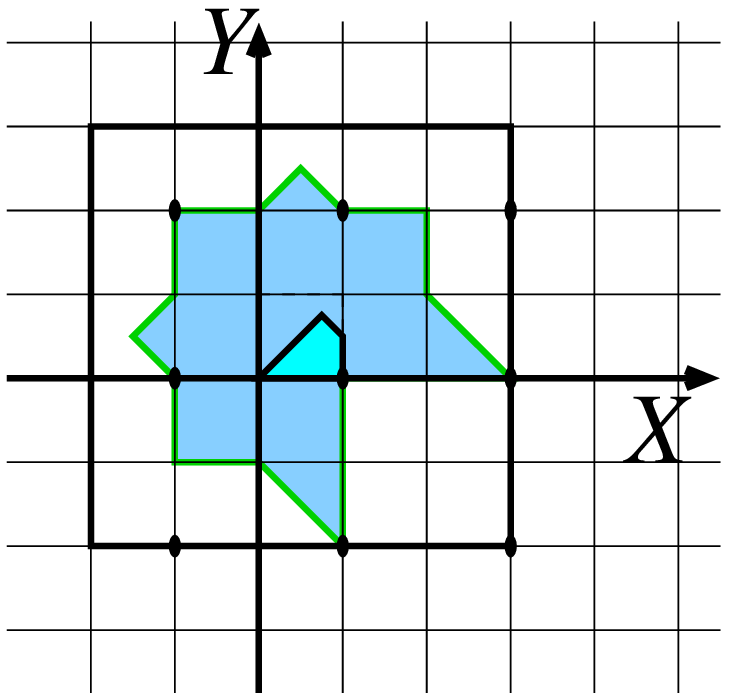}
&
\includegraphics[width=3cm]{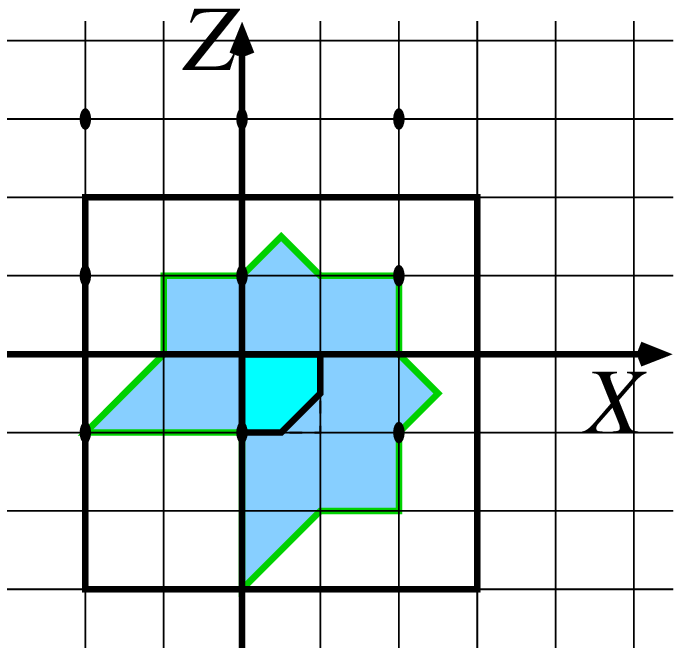} 
&
\includegraphics[width=3cm]{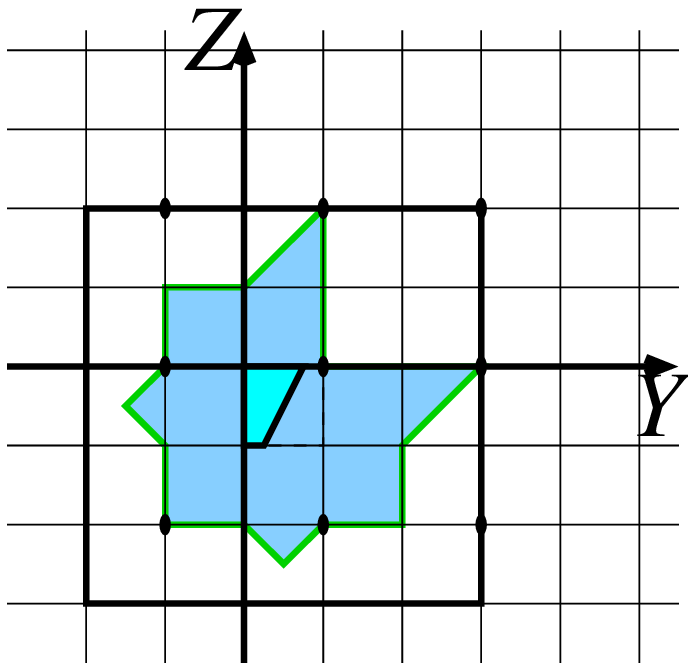}
\\
\hline
\footnotesize p2
&
\includegraphics[width=3cm]{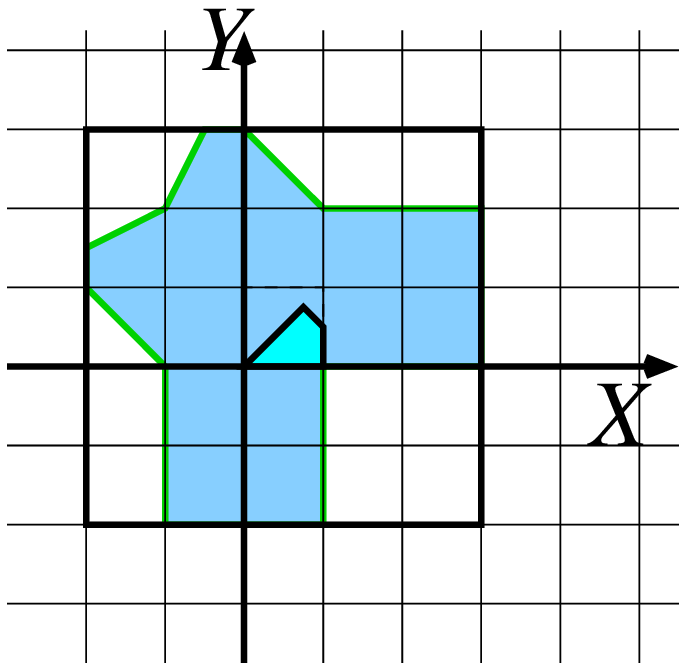}
&
\includegraphics[width=3cm]{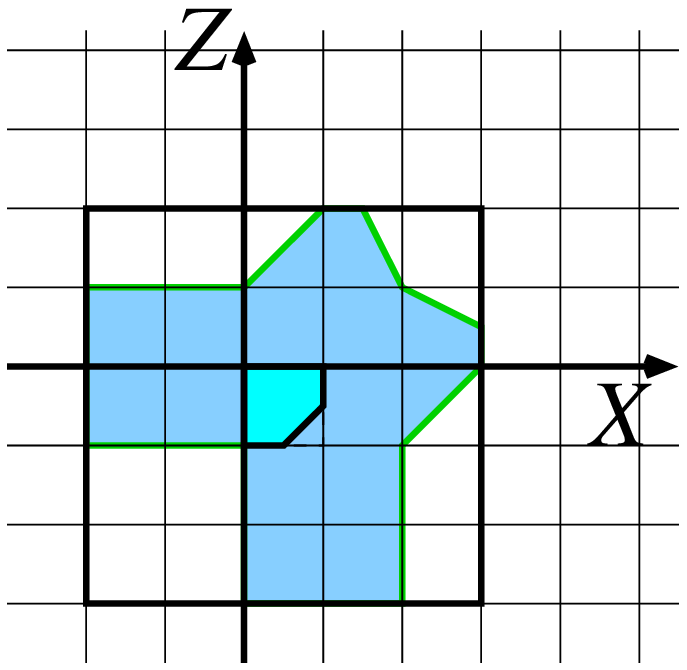}
&
\includegraphics[width=3cm]{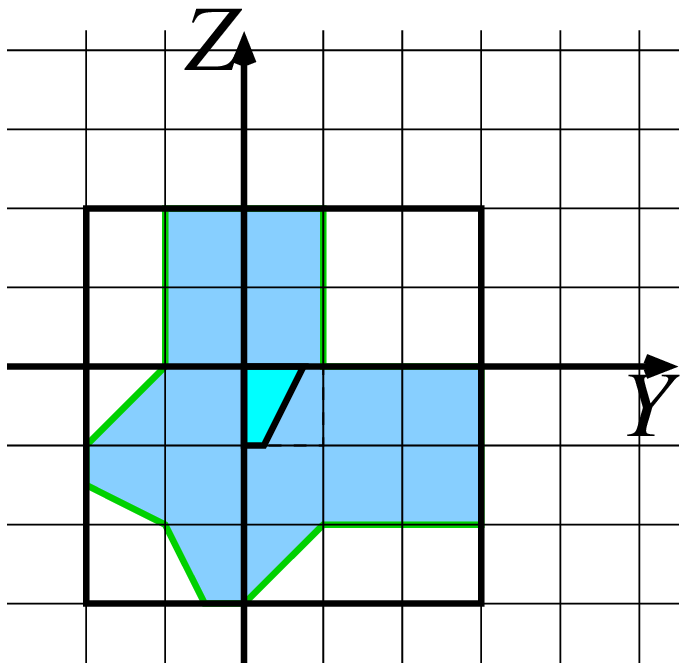}
\\
\hline
\footnotesize pg
&
\includegraphics[width=3cm]{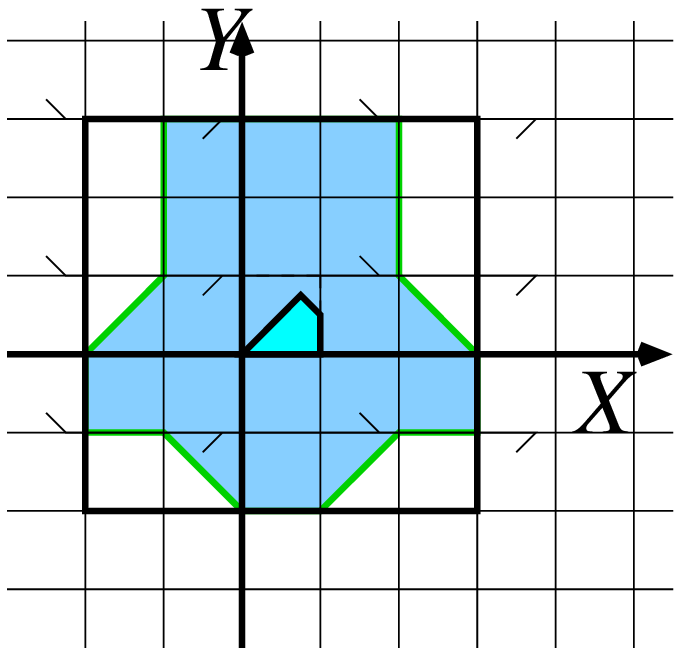}
& 
\includegraphics[width=3cm]{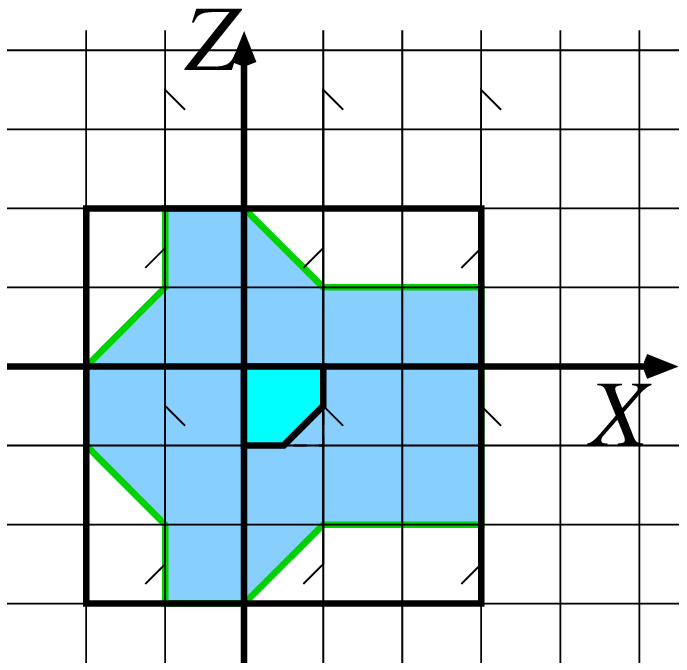}
&
\includegraphics[width=3cm]{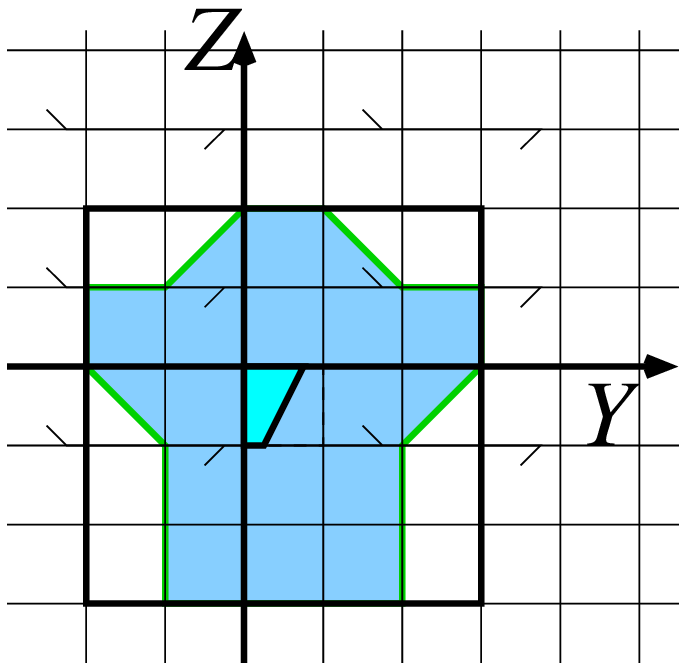}
\\
\hline
&\text{Plane XY} & \text{Plane XZ} & \text{Plane YZ} \\
\hline
\end{tabular}
\medskip
\caption{Extended Voronoi Regions of the plane groups of the quarter groups. Prototiles of types $A$ and $B$
\label{table:intersections}}
\end{table}

\begin{table}
\begin{tabular}{|c|c|c|c|}
\hline
& \multicolumn{3}{c|}{{\text{Type} {\em C}}}\\
\hline
\hline
\footnotesize pgg
&
\includegraphics[width=3cm]{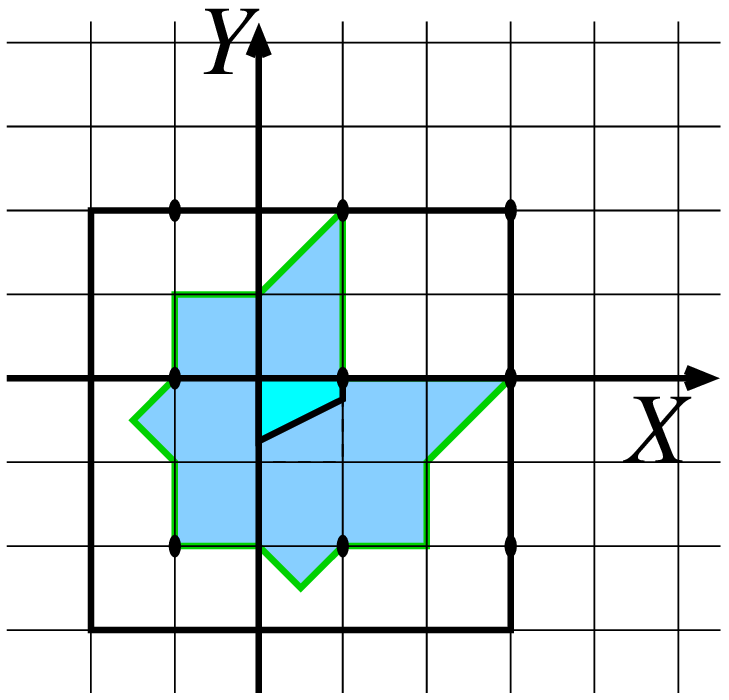}
&
\includegraphics[width=3cm]{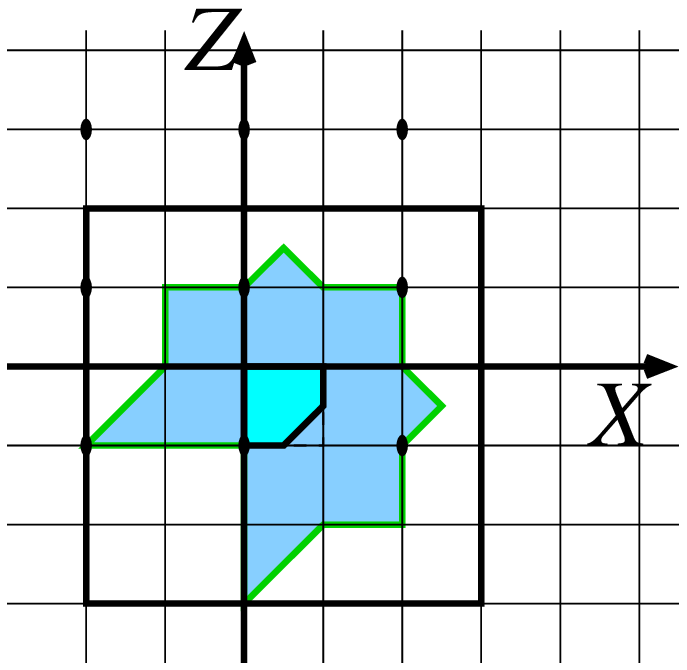}
&
\includegraphics[width=3cm]{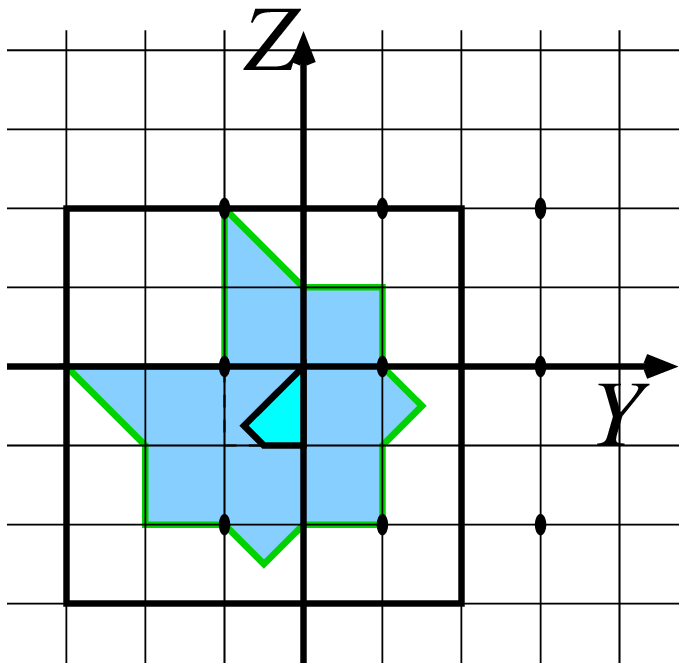}
\\
\hline
\footnotesize p2
&
\includegraphics[width=3cm]{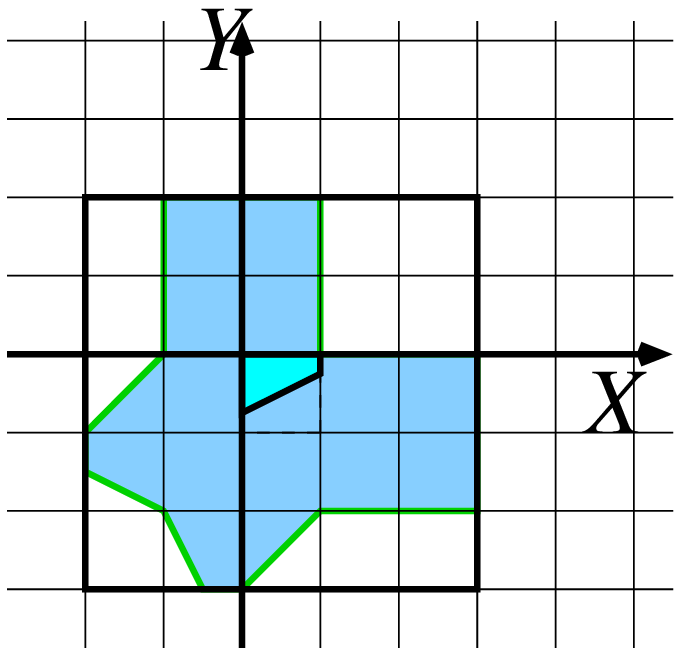}
&
\includegraphics[width=3cm]{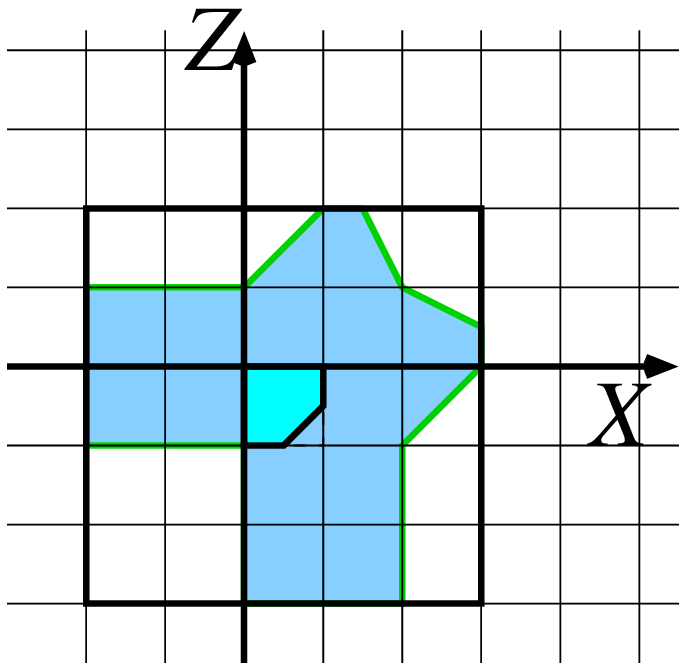} 
&
\includegraphics[width=3cm]{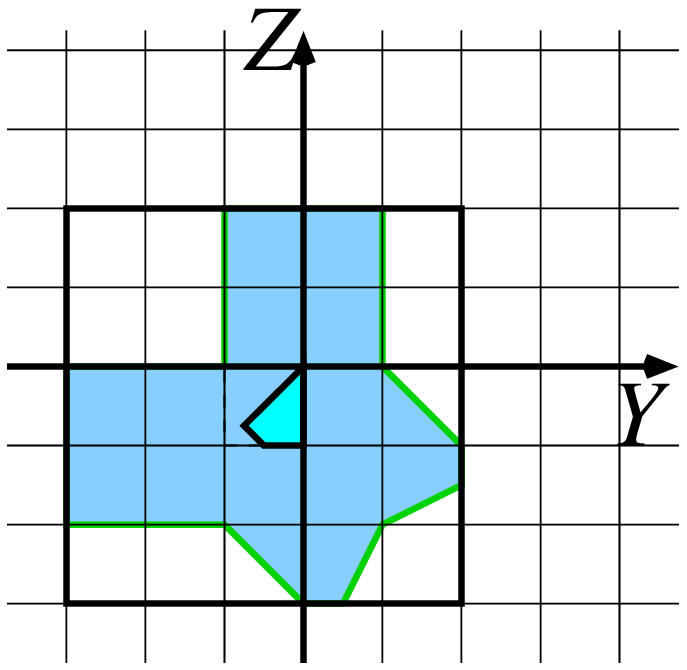} 
\\
\hline
\footnotesize pg
&
\includegraphics[width=3cm]{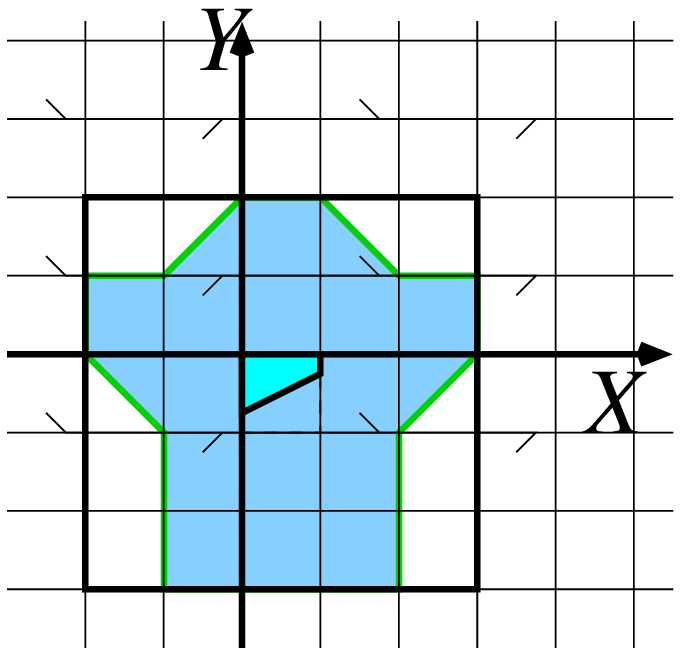}
&
\includegraphics[width=3cm]{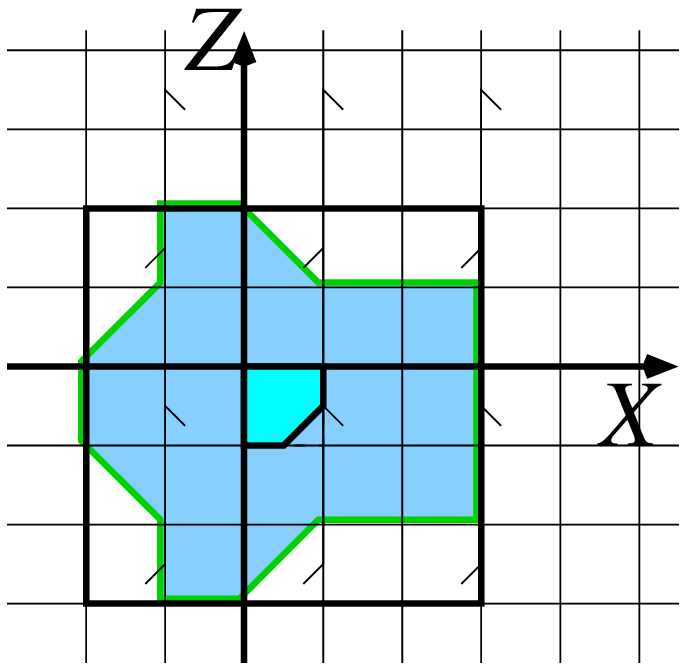}
&
\includegraphics[width=3cm]{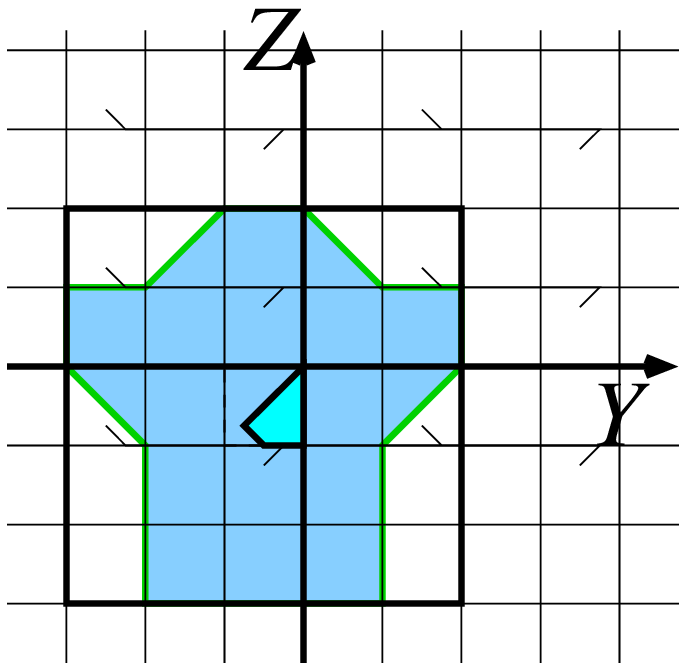}
\\
\hline
& \multicolumn{3}{c|}{{\text{Type} {\em D}}}\\
\hline
\hline
\footnotesize pgg
&
\includegraphics[width=3cm]{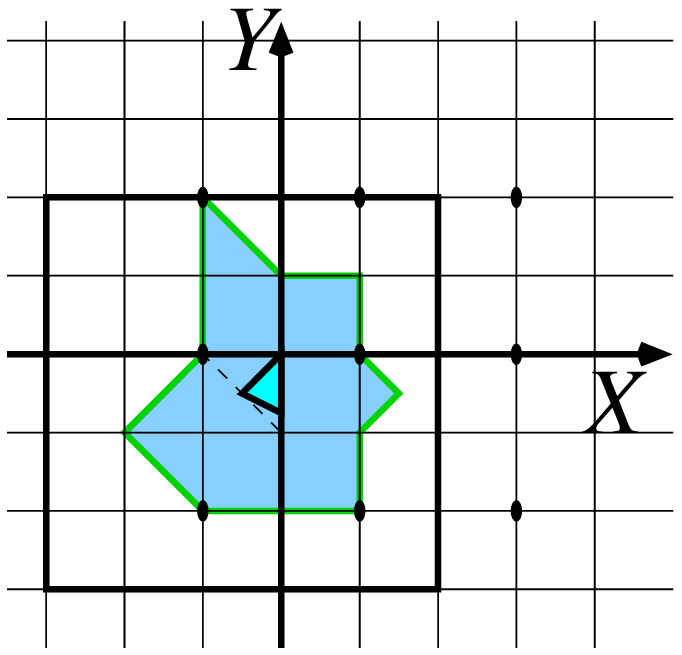}
&
\includegraphics[width=3cm]{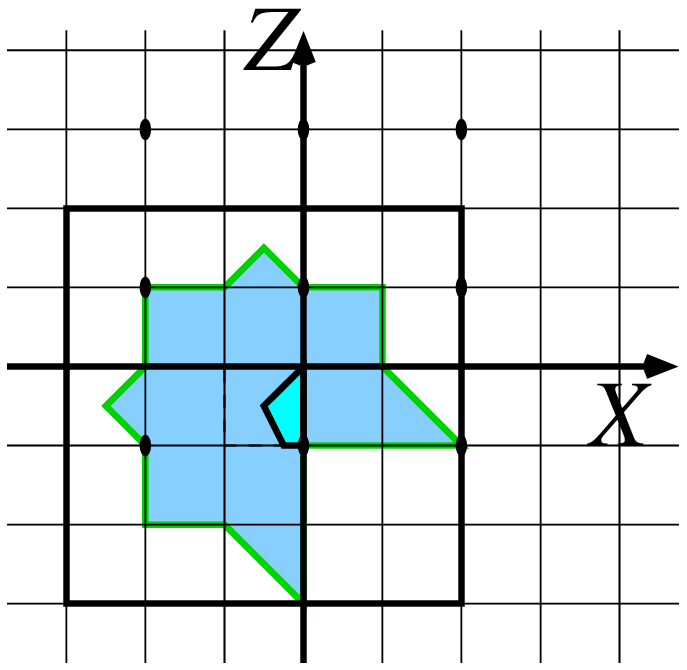} 
&
\includegraphics[width=3cm]{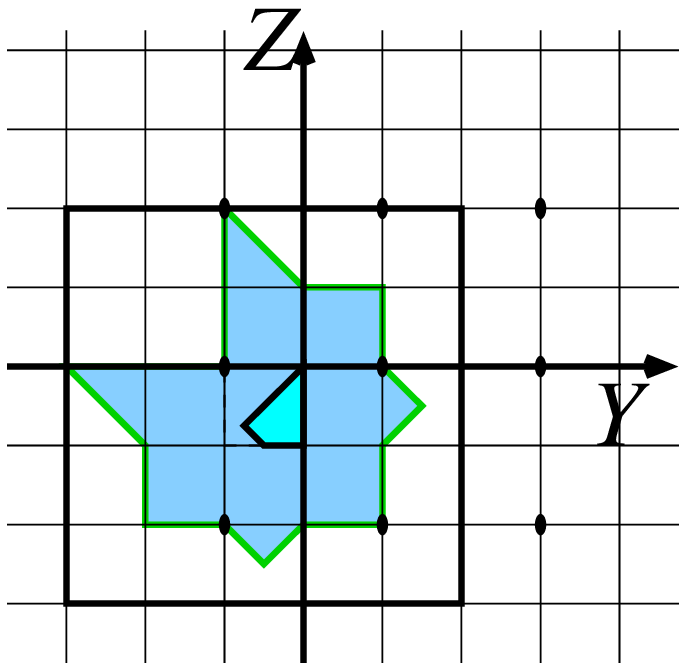}
\\
\hline
\footnotesize p2
&
\includegraphics[width=3cm]{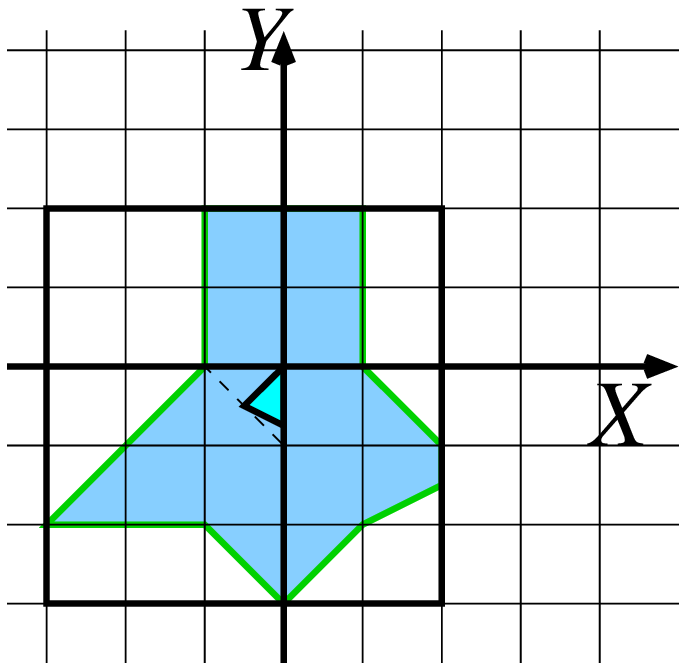}
&
\includegraphics[width=3cm]{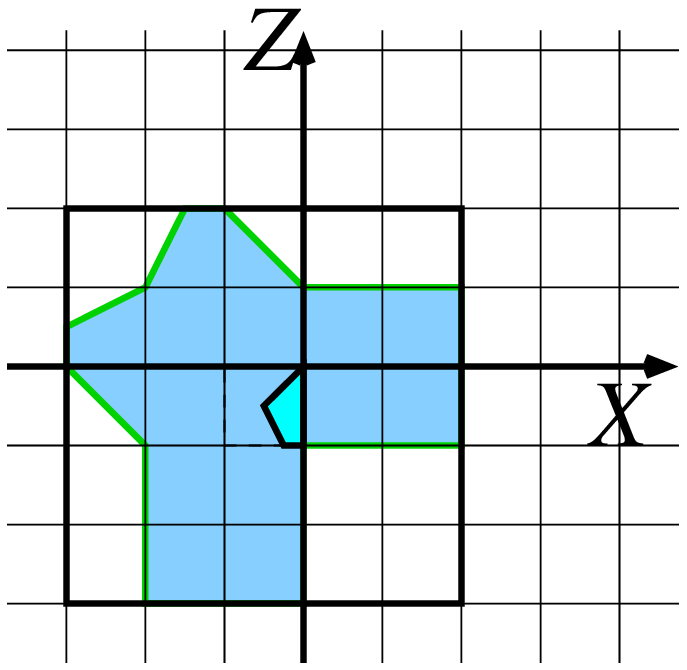}
&
\includegraphics[width=3cm]{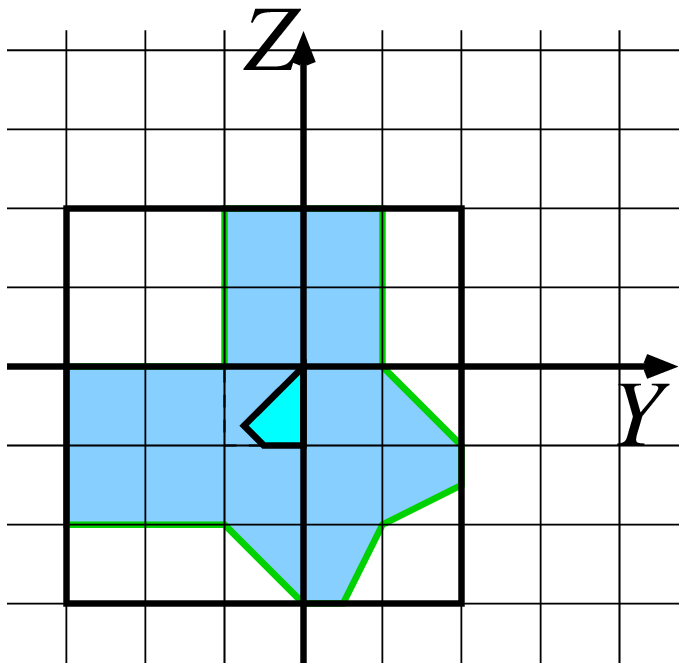}
\\
\hline
\footnotesize pg
&
\includegraphics[width=3cm]{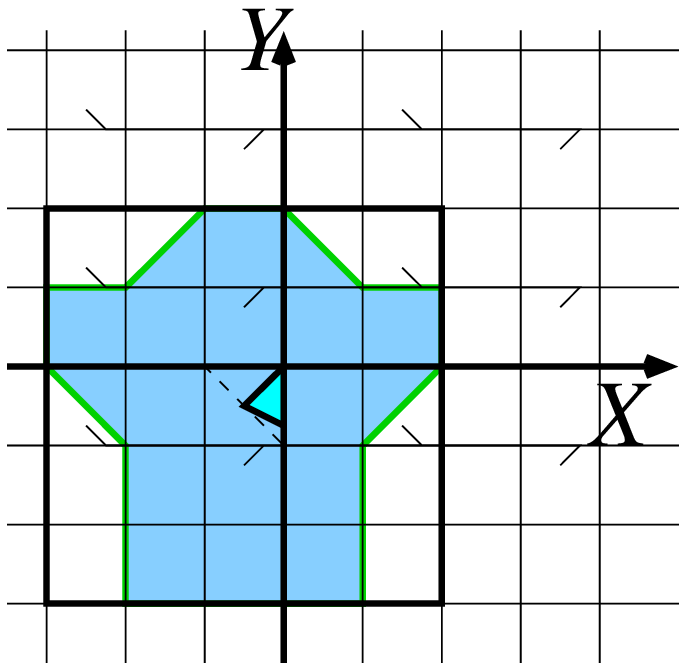}
& 
\includegraphics[width=3cm]{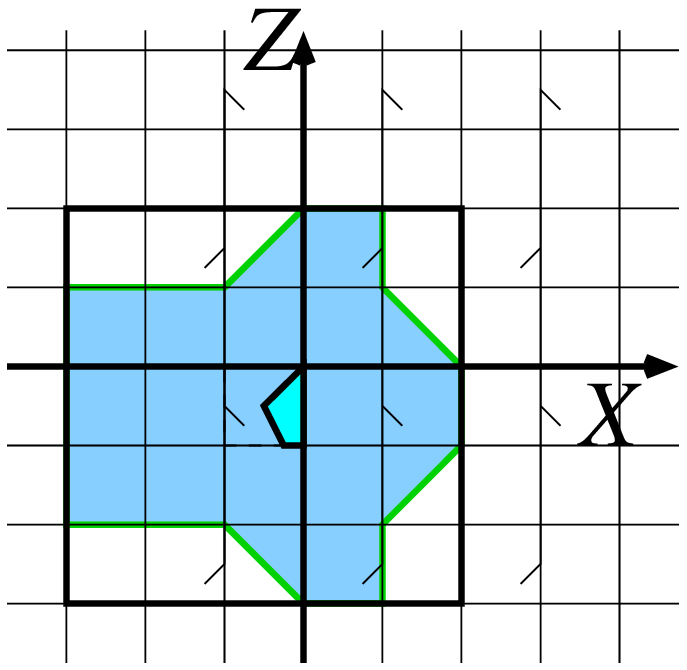}
&
\includegraphics[width=3cm]{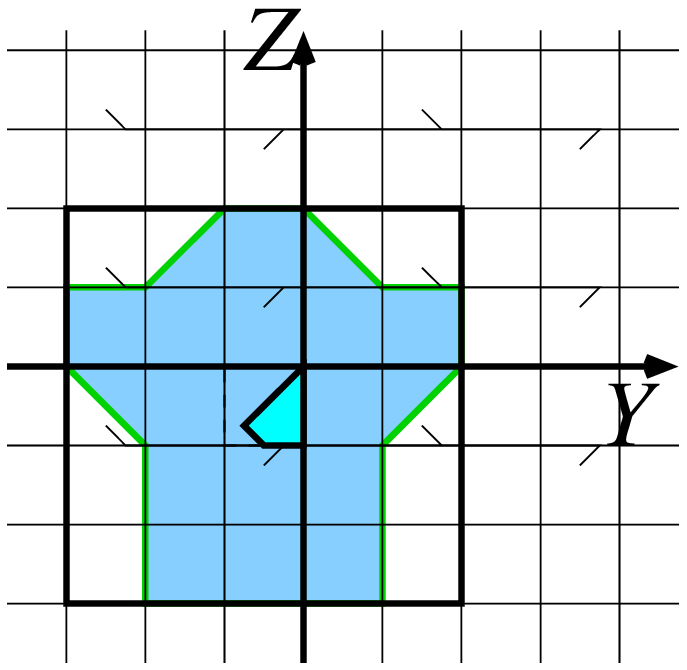}
\\\hline

& \text{Plane XY} & \text{Plane XZ} & \text{Plane YZ} \\
\hline
\end{tabular}
\medskip
\caption{Extended Voronoi Regions of the plane groups of the quarter groups. Prototiles of types $C$ and $D$
\label{table:intersections2}}
\end{table}

\subsection{Influence regions}
 \label{sec:first-bound}


By Theorem~\ref{thm:influence}, once we know  the influence region of, say, the prototile
$A_{0}$, for a group $G$, we can obtain an upper bound for the maximum number of facets of the Dirichlet stereohedra of $G$ with base point in $A_0$. We only need to count how many orbit points lie in the influence region, that is, how many tiles in the $G$-orbit of $A_0$ are listed in our influence region. (Remember that, by construction, our influence region is a union of tiles of the auxiliary tessellation $\auxtes$.)
Let us see in more detail how we compute this number.

 By definition, a tile $T$ is in $\Infl_{G}(A_{0})$ if and only if $\VorExt_{G}(A_{0})\cap \VorExt_{G}(T)\not = \emptyset$. Note that $T$ can a priori be a tile of any of the four types, $A$, $B$, $C$ or $D$. But for $T$ produce an orbit point it is necessary (but not sufficient, see below) that $T$ be of the same type, $A$, as our prototile.
Such a tile $A'$ is in $\Infl_{G}(A_{0})$ if and only if there is a tile $T'$ (of any of the four types) in $\VorExt_{G}(A_{0})\cap \VorExt_{G}(A')$. 

\begin{figure}
    \begin{center}
      \includegraphics[width=8cm]{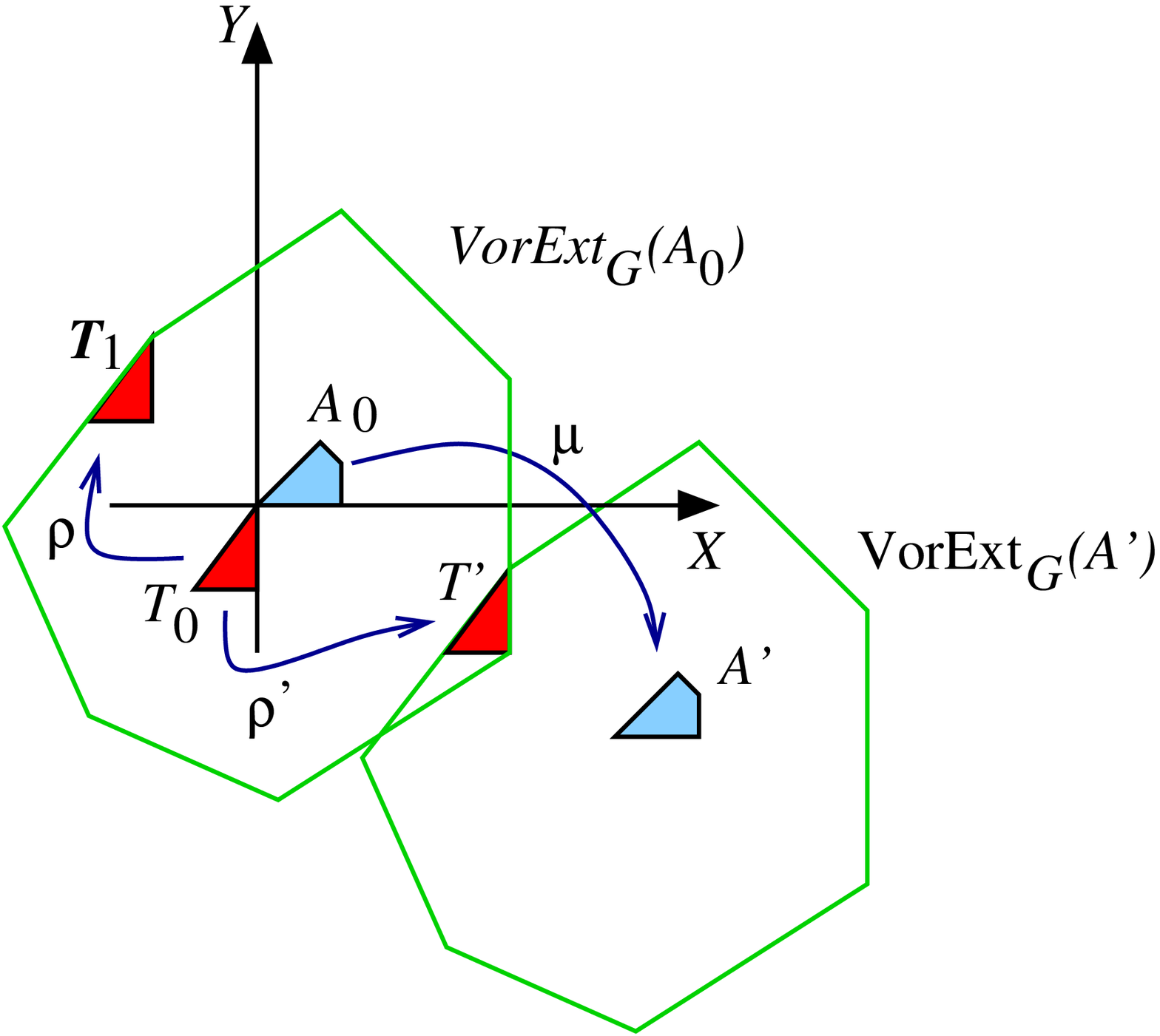}
    \end{center}
      \caption{Calculation of the influence region}
     \label{fig:influence}
\end{figure}

At this point we have to clarify one feature of our encoding of tiles that is useful here.
Remember that for us $\VorExt_{G}(A_{0})$ is a list of tiles of $\auxtes$, and that each tile $T$ of $\auxtes$ has been encoded as the transformation $\rho\in \nor(Q)$ that sends the prototile $T_0$ of the same type to $T$ (see Figure~\ref{fig:influence}). Suppose now that $A'$ is a tile of type $A$ in $\Infl_{G}(A_{0})$. That is, $A$ is  such that 
\[
T'\in\VorExt_{G}(A_{0})\cap \VorExt_{G}(A')
\]
for some $T'$. Let $\mu\in \nor(Q)$ be the transformation that sends $A_0$ to $A'$. Then, by Lemma 
\ref{lem:normalizer} (and since $\nor(Q)\subseteq\nor(G)$ for every quarter group $G$), $\VorExt_{G}(A')= \mu \VorExt_{G}(A)$. That is, there is a second tile $T_1$ in $\VorExt_{G}(A)$ such that 
$T'=\mu T_1$. Let $\rho$ and $\rho'$ denote the transformations that send the prototile $T_0$ of the appropriate type to $T'$ and $T_1$, respectively. We obviously have that (see  Figure~\ref{fig:influence} again):
\[
\mu=\rho' \circ \rho^{-1}.
\]

Therefore, all the tiles $A\subset \Infl_{G}(A_{0})$ that may possibly produce a facet in our Dirichlet stereohedra can be represented as the composition $\rho' \circ \rho^{-1}$ where $\rho$ and  $\rho'$
are two isometries (that represent tiles) in $\VorExt_{G}(A_{0})$.
 
There is a final step. Each transformation $\mu=\rho' \circ \rho^{-1}$ obtained in this way is clearly in $\nor (Q)$, but it may not be in our particular group $G$. Only if it is in $G$ it contributes one neighbor to the bound. To check this we use the coset classification of all the tiles in our initial population, shown  in Table~\ref{table:motions-groups}.
  
The influence regions for the three other  prototiles $B_0$, $C_0$ and $D_0$ are obtained in exactly the same way. Once we have the number of elements in each influence region of $G$ for each prototile, we take  the biggest of these numbers as a bound of the maximum number of facets of Dirichlet stereohedra for $G$.
 
The results of these calculations are shown in columns (1), (2), (3) and (4) of Table~\ref{table:main}. 
For each group, the rightmost number of these three columns is the maximum size of influence regions obtained. Column (1) shows the bound obtained if we neglect the diad rotations and the planar projection step. Columns (2), (3) and (4) show how the bound decreases after considering, respectively, the coordinate diad rotations, the diagonal diad rotations, and the planar projection step.

\small
\bibliographystyle{abbrv}


\end{document}